\documentclass[9pt]{article}

\usepackage{amssymb}

\usepackage{amsmath}
\usepackage{mathrsfs}
\usepackage{algorithmic}
\usepackage{algorithm}
\usepackage[pdftex]{graphicx}
\usepackage{rotating}
\usepackage{multirow}

\newtheorem{theorem}{Theorem}
\newtheorem{lemma}{Lemma}
\newtheorem{remark}{Remark}
\newtheorem{definition}{Definition}

\begin{document}

\begin{center}
  \textbf{\LARGE Purely algebraic domain decomposition methods for the incompressible Navier-Stokes equations}\vspace*{3ex}\\
  Pawan Kumar\footnote{This work was completed in part time when the
    author was provided office facilities and access to journals by
    Institut Henri Poincare, UMS 839 (CNRS/UPMC) while taking part in
    trimester program on control and PDE (Oct. - Dec, 2010),
    and Fonds de la recherche scientifique (FNRS)(Ref: 2011/V 6/5/004-IB/CS-15), Belgique. }\\
  Service de M\'etrologie Nucl\'eaire \\
  Universit\'e libre de Bruxelles \\
  Bruxelles, Belgium \\
  \textsf{pawan.kumar@u-psud.fr/kumar.lri@gmail.com}\vspace*{0.5ex}\\
\end{center}
\section*{Abstract}
In the context of non overlapping domain decomposition methods,
several algebraic approximations of the Dirichlet-to-Neumann (DtN) map
are proposed in [F. X. Roux, et. al.\newblock{Algebraic approximation
  of Dirichlet-to-Neumann maps for the equations of linear
  elasticity},\newblock{ \em Comput. Methods Appl. Mech. Engrg.}, 195,
2006, 3742-3759]. For the case of non overlapping domains,
approximation to the DtN are analogous to the approximation of the
Schur complements in the incomplete multilevel block factorization.

In this work, several original and purely algebraic (based on graph of
the matrix) domain decomposition techniques are investigated for
steady state incompressible Navier-Stokes equation defined on uniform
and stretched grid for low viscosity. Moreover, the methods proposed
are highly parallel during both setup and application phase. Spectral
and numerical analysis of the methods are also presented.
\section{Introduction}


At the core of some numerical simulations lies the problem of solving 
sparse linear systems of the form
\begin{eqnarray} \label{eqn} \mathbf{C}\mathbf{x} = \mathbf{b},
\end{eqnarray}
where $\mathbf{C} \in \mathbb{R}^{N \times N},~ \mathbf{x} \in
\mathbb{R}^{N},~ \mathbf{b} \in \mathbb{R}^{N}$.  One of the sources
of equation (\ref{eqn}) is the following time evolving Navier Stokes
equation

\begin{align}
  \frac{\partial u}{\partial t} - \nu \Delta u + (u \cdot \nabla)u + \nabla p &= f~&\text{in}~\Omega, \label{NS1} \\
  \nabla \cdot u &= 0~&\text{in}~\Omega, \label{NS2} \\
  \mathtt{B}u &= g &\text{on}~\Gamma, \label{NS3}
\end{align}
where $\nu > 0 $ is the kinematic viscosity coefficient (inversely
proportional to Reynolds number Re), $\Delta$ is the Laplace operator
in $\mathbb{R}^{d}$, $\nabla$ denotes the gradient, $\nabla \cdot $
stands for divergence, and $\mathtt{B}$ is a boundary operator. The
domain $\Omega \subset \mathbb{R}^{d}(d=2,3)$ is the bounded,
connected domain with a piecewise smooth boundary $\Gamma$. Here
$f:\Omega \rightarrow \mathbb{R}^{d}$, the boundary data given by
$g:\Gamma \rightarrow \mathbb{R}^{d}$. The system models the flow of
an incompressible Newtonian fluid such as air or water. The presence
of non-linear term $u\cdot\nabla u$ indicates presence of multiple
solutions.  We are interested in finding the velocity field $u:\Omega
\rightarrow \mathbb{R}^{d}$ and a pressure field $p:\Omega \rightarrow
\mathbb{R}$ that satisfies the equations (\ref{NS1}), (\ref{NS2}), and
(\ref{NS3}) above. Implicit time discretization such as Crank-Nikelson
scheme together with spatial discretization such as finite element
scheme \cite{elm} of the Navier-Stokes system above leads to a linear
system (\ref{eqn}) where the matrix $\mathbf{C}$ has the following
block $2 \times 2$ partitioned form

$$
\mathbf{C}=\left[\begin{array}{c|c} D & E^{T} \\ \hline -E & 0
  \end{array}
\right].$$ Here $E^{T}$ is the discrete gradient and $-E$ is the
negative divergence operator. The structure of the matrix $D$ depends
on the nonlinear algorithm; for Picard iteration $D$ is block diagonal
with each diagonal block being the discrete convection-diffusion
operator. For Newton iteration, $D$ has more complex structure
\cite{elm}. The matrix $\mathbf{C}$ is indefinite and non-symmetric.

In a real life simulation, the coefficient matrix $\mathbf{C}$ is
usually large and sparse, and the usual direct methods \cite{dav} are
costly both in terms of CPU time and storage requirements. A common
preference is to use Krylov subspace methods with some suitable
preconditioning. A preconditioner $\mathbf{B}$ is an approximation to
the coefficient matrix $\mathbf{C}$ such that the preconditioned
operator $\mathbf{B^{-1}C}$ has ``favourable'' spectrum that is
essential for a fast convergence of Krylov subspace based iterative
methods \cite{saad96}. In general, the number of iterations required
for convergence is less when the eigenvalues are clustered towards one
and they are away from zero. On the other hand, the time required to
setup the preconditioner, and the cost of applying it during
the iteration phase should not be too demanding.

Most of the classical and recent preconditioners for the Navier-Stokes
systems are approximate block factorization (ABF). Classical
pressure correction methods are SIMPLE (Semi-Implicit Method for
Pressure Linked Equations), and their modifications SIMPLEC, and
SIMPLER \cite{pat2,per,per2,pat}. A promising class of method based on
approximation of the Schur complement ($S = ED^{-1}E^{T}$) is a
pressure convection diffusion (PCD) preconditioner, where, the Schur
complement is approximated by
\begin{eqnarray}
  S\approx \widehat{S} = EE^{T}D_{p}^{-1}.
\end{eqnarray}
Thus, the approximation $\widehat{S}$ is obtained by first
approximating $D$ by $D_{p}$ and then commutating $D_{p}^{-1}E^{T}$ to
get $D_{p}$ on the right. Here $D_{p}$ is the discrete
convection-diffusion operator projected on the pressure space
\cite{elm}. The method leads to convergence rates that are independent
of the mesh size but deteriorates with Reynolds number higher that 100
as confirmed in the numerical experiments section
\cite{elm99,elm02,kay02}. A modification of PCD is the least squares
commutator (LSC) preconditioner, where the construction of the
discrete convection-diffusion operator projected on the pressure space
is automated by solving the normal equations associated with a least
square problem derived from the commutator in \cite{elm68}. Another
approach is based on the Hermitian or Skew Hermitian Splitting (HSS)
and Dimensional Splitting (DS) of the problem along the components of
the velocity field and its relaxed version are introduced in
\cite{ben2,ben4}. Here, HSS has not been implemented efficiently for
Oseen problems and DS preconditioning suffers poor convergence on low
viscosity problems on stretched grid.  In general, these methods
belong to class of preconditioners where the preconditioner has block
$2\times2$ form, and for parallelism, they rely on the scalability of
the inner solvers for (1,1) and Schur complement block in the
preconditioner.

In this paper, we concern ourselves with substructured domain
decomposition (DD) based preconditioner. On contrary to overlapping DD
methods, the non-overlapping DD methods tries to approximate the
so-called Dirichlet to Newmann (DtN) map
\cite{Lions,nataf98,nataf1994,nat,gan}. For its algebraic counterpart,
approximation of DtN map is related to the approximation of the Schur
complement \cite{nat,rou1,rou2}. In the non-overlapping DD method
considered in this paper, the required substructuring (partitioning)
is obtained by a graph partitioner \cite{kar3,saad2003} thus leading
to a purely algebraic method where the graph of the matrix is
sufficient and nothing else is assumed of the computational domain and
boundary conditions. The partitioner finds a set of nodes (separator),
removal of which leaves the graph disconnected into as many
disconnected components as required. That is, the matrix $\mathbf{C}$
above is transformed to $\mathbf{P}^{T}\mathbf{C}\mathbf{P}$ where
$\mathbf{P}$ is the permutation matrix that resuffles the rows and
columns.  Thus, the resulting permuted matrix can now be partitioned
into block $2 \times 2$ form such that the (1,1) block is block
diagonal with as many blocks as desired. We notice here that such
reordering techniques are popular and almost always taken into
consideration when degisning direct \cite{dav} or hybrid methods
\cite{saad2003,ali} to reduce the amount of fill-in and to enhance
parallelism \cite{saad2003}. The main contribution of our work is the
degisn of a parallel computation of Schur complement.  In \cite{rou1},
some algebraic approximations are considered for the Schur complement;
in one of the methods the global Schur complement is approximated in
patches. The computation of these patch Schur complements ultimately
leads to an approximation of the global Schur complement. In this
paper, we propose some modifications in the approximation of the Schur
complements. Rather than building a patch around a node, we consider
an aggregated set of nodes, and build a mini Schur complement
approximation (MSC) for all the nodes of the aggregate at once. Apart
from this basic modification, we propose to construct ``patches''
based on the numbering of the nodes rather than on the edge
connections. In other words, patches consists of closely numbered
nodes rather than closely connected ones. This new approach leads to
much faster approximations. The method is purely algebraic and takes
matrix and right hand side as an input, and easily integrated in a
non-linear solver. Compared to two state-of-the-art ABF methods namely
PCD and LSC, the proposed methods are attractive for several reasons:
\begin{itemize}
\item The setup and application phase of the preconditioners are
  massively parallel
\item The new methods converges faster and in particular, compared to
  PCD and LSC methods, they perform significantly better for Reynolds
  number larger than 100.
\item The diagonal blocks of the (1,1) blocks are approximated by
  incomplete LU factorizations leading to sparse factors for the
  preconditioner.
\end{itemize}
Although, the method can be tried on any problem, we concern ourselves
with problems steming from Navier-Stokes equation for Reynolds number
ranging from 10 to 3000.

The remainder of this paper is organized as follows. In section (2),
we explain briefly the PCD and LSC methods.  In section (3), we
introduce the new methods based on mini Schur complements, the
importance of overlapping the patches will be studied. In section (4),
we discuss the parallelism and implementation aspects for the new
methods. Finally, in section (5), we present the numerical experiments
and we compare, the fill factor, the iteration count, and the
robustness of the methods.

\section{Some preconditioners for the incompressible Navier-Stokes equations}
In this section, we briefly describe the pressure convection diffusion
(PCD) and least squares commutator (LSC) preconditioners. These methods 
will serve as a benchmark methods for comparison.

\subsection{Pressure Convection Diffusion}
Let the discrete convection diffusion operator on the velocity space 
be defined as follows
\begin{eqnarray}
  \mathscr{L} = -\nu \partial^{2}+w_{h}\cdot \nabla
\end{eqnarray}
where $w_{h}$ is the approximation to the discrete velocity in the
recent Picard iteration. Consider a similar operator $\mathscr{L}_{p}$
is defined on the pressure space
\begin{eqnarray}
  \mathscr{L}_{p} = (-\nu \partial^{2}+w_{h}\cdot \nabla)_{p}
\end{eqnarray}
Let $\epsilon$ denote the commutator of both these operator with the
gradient operator as follows
\begin{eqnarray} \label{eqn_com} \epsilon = \mathscr{L} \nabla -
  \nabla \mathscr{L}_{p}
\end{eqnarray}
The discrete analog of the commutator is given as follows
\begin{eqnarray}
  \epsilon_{h} = (Q^{-1}D)(Q^{-1}E^{T}) - (Q^{-1}E^{T})({Q_{p}}^{-1}D_{p}).
\end{eqnarray}
Here $Q$ and $Q_{p}$ are velocity and pressure mass respectively. The
transformation from integrated to nodal values is done by the inverse
operations of velocity and pressure mass matices. Here $D_{p}$ is the
discrete convection-diffusion operator on pressure space. Following
approximation to the Schur complement
\begin{eqnarray} \label{pcd}
  E D^{-1} E^{T} \approx  E Q^{-1} E^{T} {D^{-1}}_{p} Q_{p}
\end{eqnarray}
is obtained by assuming that $(ED^{-1}Q)\epsilon({F_{p}}^{-1}Q_{p})$
is close to zero. We observe that $E$ being discrete divergence oeprator and
$E^{T}$ being the gradient operator, $EE^{T}$ is the discrete Laplacian
and $EQ^{-1}E^{T}$ is the scaled discrete Laplacian. Here
$EQ^{-1}E^{T}$ being expensive is replace by its spectrally equivalent
pressure mass matrix $A_{p}$, thus leading to a final approximation
$\widehat{S}$ of the Schur complement to be
\begin{eqnarray}
  \widehat{S} = -A_{p}D_{p}^{-1}Q_{p}
\end{eqnarray}
The PCD preconditioner denoted by $\mathbf{B_{PCD}}$ is defined as follows
$$
\mathbf{B_{PCD}}=\left[\begin{array}{c|c} D & E^{T} \\
    \hline 0 & \widehat{S}
  \end{array}
\right].
$$

\subsection{Least Square Commutator}
One of the drawbacks of the PCD method is that we need to construct the
convection-diffusion operator $D_{p}$ projected on the pressure space
which essentially requires  an full understanding of the underlying
discretization scheme and other impelmentation details. In \cite{elm99,elm68}, 
Elman et. al. propose to find $D_{p}$ automatically by solving a least square problem 
of the form min$\|\epsilon_{h}\|_{Q}$, i.e., by minimizing
\begin{eqnarray}
  min\|(Q^{-1}DQ^{-1}E^{T})_{j} - Q^{-1}E^{T}Q_{p}^{-1}(D_{p})_{j}\|_{Q},
\end{eqnarray}
where $\|\cdot\|_{Q}$ is $\sqrt{x^{T}Qx}$ norm, and $(K)_{j}$ for any matrix $K$ 
denotes the $j^{th}$ column of $K$. The normal equations
associated with this problem are given as follows
\begin{eqnarray}
  Q_{p}^{-1}EQ^{-1}E^{T}Q_{p}^{-1}(D_{p})_{j} = (Q_{p}EQ^{-1}DQ^{-1}E^{T})_{j}
\end{eqnarray}
which gives
\begin{eqnarray*}
  D_{p} = Q_{p}(EQ^{-1}E^{T})^{-1}(EQ^{-1}DQ^{-1}E^{T}).
\end{eqnarray*}
Substuting the expression for $D_{p}$ in (\ref{pcd}), we obtain an
approximation to the Schur complement as follows
\begin{eqnarray*}
  ED^{-1}E^{T} \approx \widehat{S}= (EQ^{-1}E^{T})(EQ^{-1}DQ^{-1}E^{T})^{-1}(EQ^{-1}E^{T}).
\end{eqnarray*}
Solving with this approximation requires two discrete poisson solve
(scaled Laplacian) which can be handled efficiently by multigrid
methods \cite{rug,stu,tro,bra}, and in contrast to other Schur
complement based methods, we only need a matrix vector product with
$D$. The LSC preconditioner denoted by $\mathbf{B_{LSC}}$ is
defined as follows
$$
\mathbf{B_{LSC}}=\left[\begin{array}{c|c} D & E^{T} \\
    \hline 0 & \widehat{S}
  \end{array}
\right].
$$

\section{Mini Schur complements}
\subsection{Graphical view} \label{graphview}

In this section, we introduce an aggregation based mini Schur
complement.  Consider again the following block $2\times2$ partitioned
matrix

\[ \label{twoBytwo} \mathbf{C}=\left[
  \begin{array}{c|c}
    D & E \\ \hline
    F & G
  \end{array}\right]
\]

Let $S = G - FD^{-1}E$ denote the global Schur
complement. Furthermore, let $\mathscr{D}$, $\mathscr{E}$,
$\mathscr{F}$, and $\mathscr{G}$ denote the set of vertices
corresponding to the adjacency graph of matrices $D$, $E$, $F$, and
$G$ respectively. Also, assume a local numbering of nodes in
$\mathscr{D}$ and $\mathscr{G}$, and for simplicity, we assume that
the number of nodes in $\mathscr{D}$ is greater than the number of
nodes in $\mathscr{G}$. The MSCs are constructed in the following
steps.
\begin{enumerate}
\item Choose a set of aggregated nodes in graph $\mathscr{G}$,
  $\mathscr{V} =
  \{\mathscr{G}_{p_{1}},\mathscr{G}_{p_{2}},\dots,\mathscr{G}_{p_{k}}
  \}$, $\mathscr{G}_{p_{i}}\subset \mathscr{G}$, $\mathscr{G}_{p_{i}}
  \cap \mathscr{G}_{p_{j}}=\emptyset$ for $i \neq j$, $\mid
  \mathscr{G}_{p_{i}}\mid=p_{i}$ and
  $\cup_{i}\mathscr{G}_{p_{i}}=\mathscr{G}$. One possible choice of
  aggregation is simply choosing the nodes with consecutive numbering
  as follows
  \begin{eqnarray} \label{aggrV} \mathscr{V} = \{ \{1,2,\dots,p_{1}
    \}, \{p_{1}+1, p_{1}+2, \dots, p_{1}+p_{2}\}, \dots \}.
  \end{eqnarray}
  Note here that the aggregation $\mathscr{V}$ is done based
  on the ``numbering'' of the nodes in the grid, rather than on the
  ``closeness'' of the nodes determined by looking at the edge
  connectivity between the nodes.  For example, in the case of 2D $n
  \times n$ grid, the node numbered $i$ could be at a distance (or
  path length) $i+n$ from node numbered $i+1$.

\begin{remark}
  We notice here that a simple generalization of the above aggregation
  scheme is obtained when we allow overlapping between the aggregates,
  i.e., for \textbf{overlapped aggregation} scheme we consider
  $\mathscr{V} =
  \{\mathscr{G}_{p_{1}},\mathscr{G}_{p_{2}},\dots,\mathscr{G}_{p_{k}}
  \}$, $\mathscr{G}_{p_{i}}\subset \mathscr{G}$, $\mathscr{G}_{p_{i}}
  \cap \mathscr{G}_{p_{j}} \neq \emptyset$ for $i \neq j$, $\mid
  \mathscr{G}_{p_{i}}\mid=p_{i}$ and
  $\cup_{i}\mathscr{G}_{p_{i}}=\mathscr{G}$.  As in the case of
  overlapped Schwarz methods, overlapping increases sharing of
  information between the MSCs thus leading to an improved
  approximation of the global Schur complement. This overlapped
  aggregation scheme could be applied to all the methods that follow
  thus leading to an improvement in the approximation of the method
  concerned. However, it is to be noted that overlapping leads to lack
  of parallelism during the solve phase since the global approximated
  Schur complement is no longer block diagonal.
\end{remark}

\item Next, we have following three possible approximation schemes
  \begin{enumerate}
  \item \textbf{Mini Schur complements based on edge connectivity
      (MSCE):} For each aggregated nodes $\mathscr{G}_{p_{i}}$,
    identify a set of nodes $\mathscr{D}_{r_{i}}$ in the set
    $\mathscr{D}$ such that for each node of $\mathscr{G}_{p_{i}}$,
    there exist a node in $\mathscr{D}_{r_i}$ within a path length of
    $r_i$. That is, we can reach a node in $\mathscr{G}_{p_i}$ from at
    least one node in $\mathscr{D}_{r_i}$ by traversing a path of
    length less than or equal to $r_i$ in the adjacency graph of
    matrix $\mathbf{C}$.  The edges between the two aggregates
    $\mathscr{D}_{p_{i}}$ and $\mathscr{G}_{p_{i}}$ are denoted by
    $E_{i}$ (incoming edge) and $F_{i}$ (outgoing edge). When this
    method is defined for the overlapping aggregates, we shall call it
    OMSCE (overlapped mini Schur complement based on edge
    connectivity).
  \item \textbf{Mini Schur complements based on numbering (MSCN):} For
    each aggregated nodes $\mathscr{G}_{p_{i}}$, identify a set of
    nodes in the graph $\mathscr{D}$ by setting $\mathscr{D}_{p_{i}}$
    = $\mathscr{G}_{p_{i}}$. Remember that the graphs
    $\mathscr{D}_{p_{i}}$ and $\mathscr{G}_{p_{i}}$ have local
    numbering of the nodes. Thus, we identify aggregates which have
    same numbering in their respective graph. When this method is
    defined for the overlapping aggregates, we shall call it OMSCN
    (overlapped mini Schur complement based on numbering).
  \item \textbf{Lumped approximation of Schur complement (Lump):} Do
    not do anything for Lump method. When this method is defined for
    the overlapping aggregates, we shall call it OLUM (overlapped
    Lumped method.)

  \end{enumerate}

\item Here again we have three cases as follows
  \begin{enumerate}
  \item \textbf{Computation of mini Schur complements for MSCE:}
    Assemble the nodes of the aggregate $\mathscr{G}_{p_i}$ and the
    edge connections between the nodes in the matrix
    ${G}_{p_i}$. Similarly, Assemble the nodes of aggregates
    $\mathscr{D}_{r_i}$ and the edge connections between the nodes in
    the matrix ${D}_{r_i}$. For a non-symmetric matrix $\mathbf{C}$,
    the corresponding graph is considered as a directed graph where
    the entries below the diagonal of the matrix may represent the
    incoming edge and those above the diagonal are outgoing edges. Let
    ${F}_{i}$ be the matrix which stores the incoming edges and
    ${E}_{i}$ denotes the outgoing edges. We consider the following
    sub matrix for the $ith$ aggregate
    
    \[
    \mathbf{C_{p_i,r_i}}=\left[
      \begin{array}{c|c}
        {D}_{r_i}  & {E}_{i} \\ \hline
        {F}_{i} &  {G}_{p_i}
      \end{array}\right].
    \]
    
    The corresponding $i^{th}$ mini Schur complement is given as
    follows
    \begin{eqnarray*}
      {S}_{i} = {G}_{p_i} - {F}_{i}({D}_{r_i})^{-1}{E}_{i}
    \end{eqnarray*}
    Here $S_{i}$ is called the mini Schur complement. In case, the
    matrix $D_{r_i}$ is large we may approximate $S_{i}$ by
    \begin{eqnarray}
      {S}_{i} = {G}_{p_i} - {F}_{i}({D}_{r_i})^{-1}{E}_{i}\mathbf{1}
    \end{eqnarray}
    where $\mathbf{1}=[1,1,1,\dots,1]^{T}$. We call this method MSCER,
    where R stands the rowsum.
  \item \textbf{For MSCN: }

    We have shown the computation of Schur complement for MSCE but we
    can proceed in a similar way to compute the mini Schur complement
    for MSCN method by replacing ${D}_{r_i}$ by ${D}_{p_i}$ above. In
    the sections that follows, we shall give a matrix view of the MSCN
    method. A simple generalization of MSCN method is to take the
    ${D}_{r_i}$ slightly larger and not necessarily the same size as
    ${G}_{r_i}$. A colsum approximation of the Schur complement as
    done above for MSCE method will be called MSCNR.
    
  \item \textbf{For Lump: } The corresponding $i^{th}$ mini Schur
    complement for the lumped approximation is given as follows
    \begin{eqnarray*}
      {S}_{i} = {G}_{p_i}
    \end{eqnarray*}
  \end{enumerate}
  
\item Extract all the columns corresponding to the nodes of the
  aggregate ${S}_{i}$ corresponding to the nodes in the aggregate
  $\mathscr{G}_{p_i}$ and put them in the corresponding columns in the
  global Schur complement. This step is same for MSCN and Lump method.
  In other words, set $S$ = blkDiag($S_{i}$).
\end{enumerate}

\subsection{Illustration of MSCN method with an example}
In this section, we illustrate the MSCN method graphically for an
small example, but for implementation, we refer the reader to
Algorithms (\ref{pseudocode_build}) and (\ref{pseudocode_solve}).

\begin{figure}[h]

  \caption{ \label{aggregation} An example of an
    aggregation. Left: An
 aggregated set of $\mathscr{G}$ containing
    two nodes indicated in bold spheres. Right: An aggregated set of
    $\mathscr{D}$ containing two nodes indicated in spheres with
    pattern. }
 
  \begin{center}
    $ \begin{array}{cc}
      \includegraphics[scale=0.3]{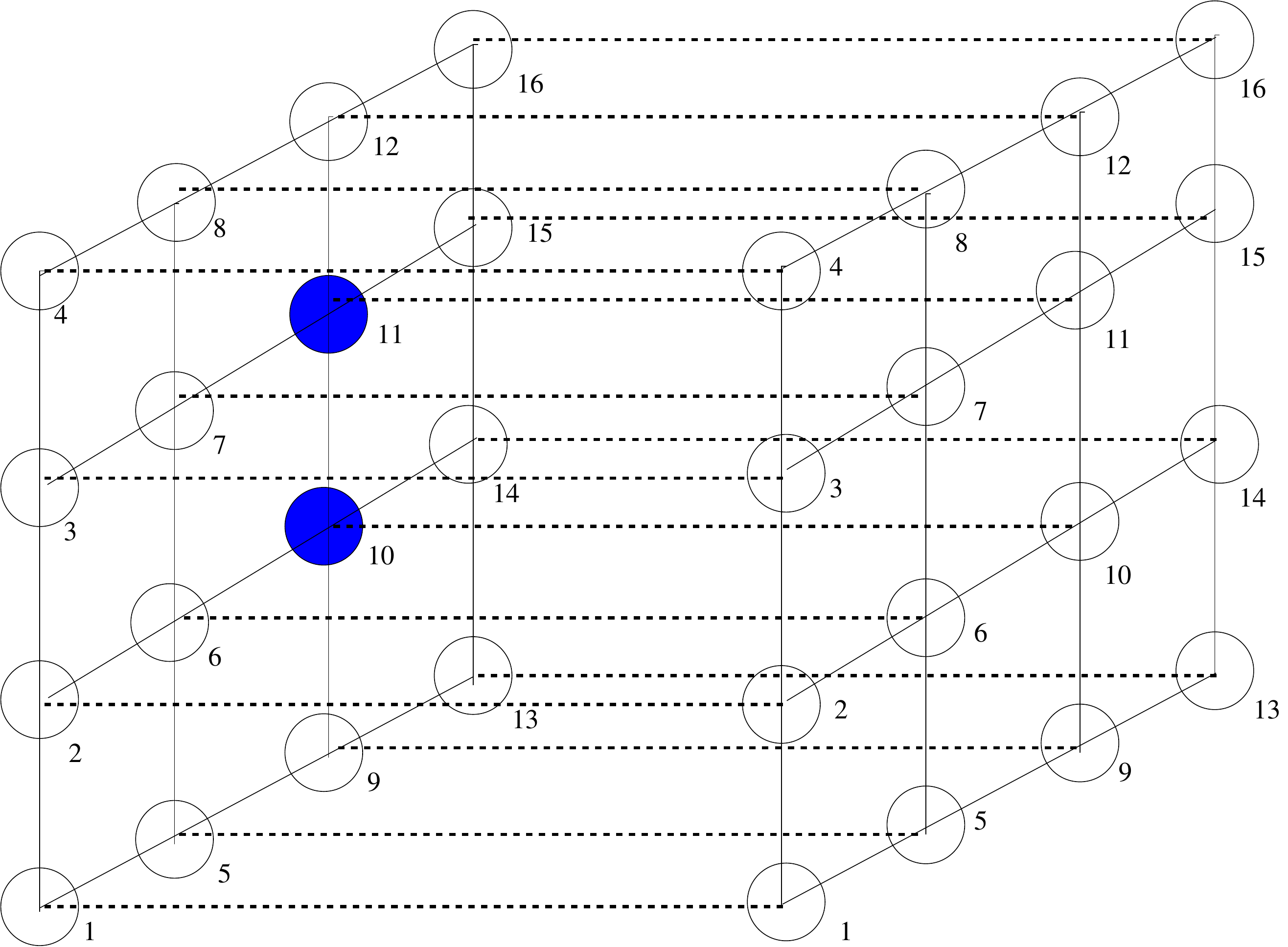} & 
      \includegraphics[scale=0.3]{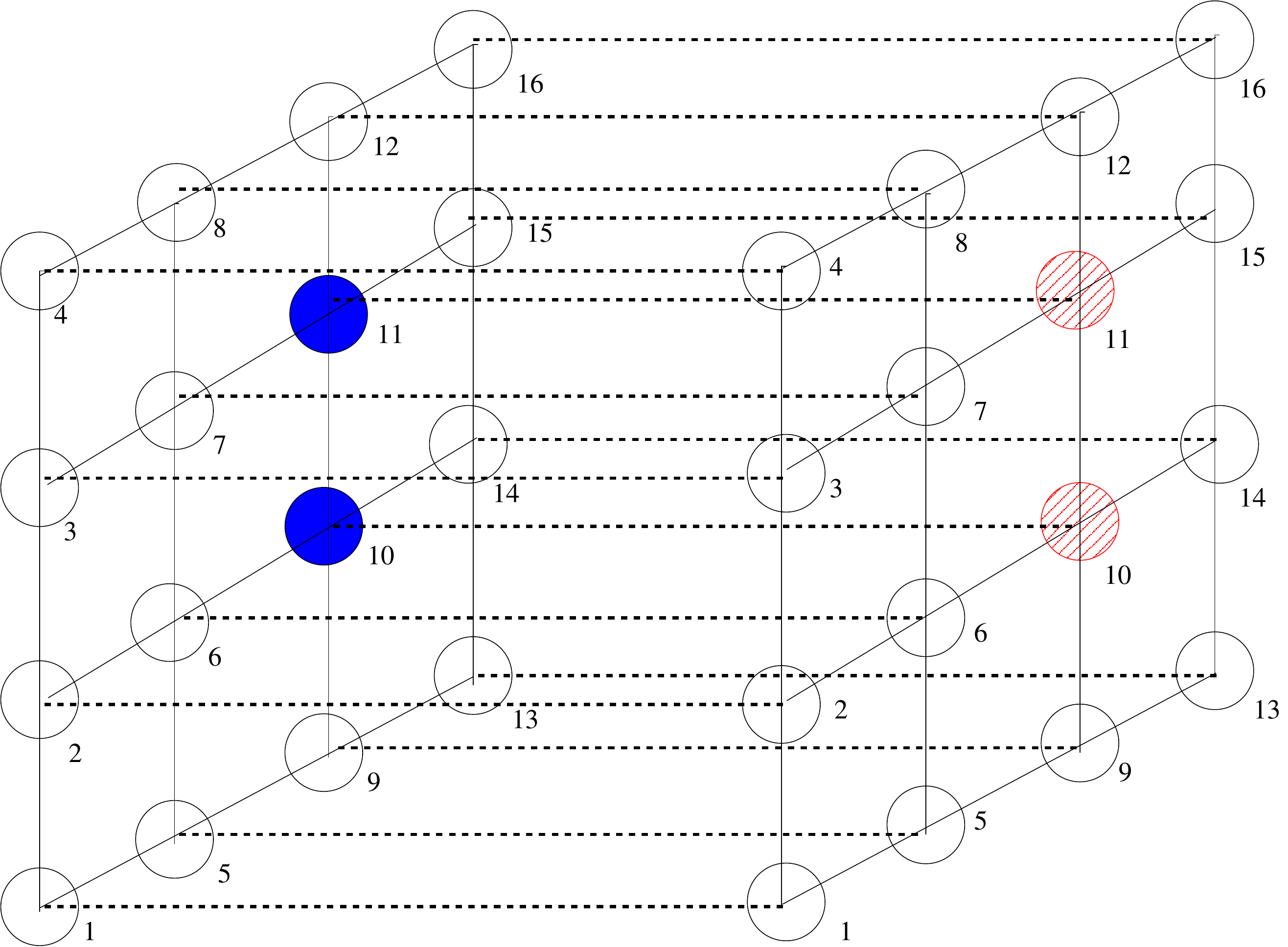} \\
    \end{array}$
  \end{center}
 
\end{figure}

In Figure (\ref{aggregation}), the two sub-figures illustrate the
steps involved in forming a mini-Schur complement.  In each sub
figure, there are two vertical planar graphs, and each graph keeps
their own local numbering of the nodes. The left vertical planar graph
within each sub figure corresponds to matrix $G$, and the right graph
corresponds to the matrix $D$. The whole graph corresponds to the
matrix $\mathbf{C}$. The steps involved in building a mini-Schur
complement, in this case, are illustrated as follows

\begin{enumerate}
\item Choose aggregates a priori for graph of G. One possible aggregate is
  \begin{eqnarray*}
    \mathscr{V} = \{ \{1, 2\}, \{3, 4, 5\}, \{6, 7, 8, 9\}, \{10, 11\},
    \{12, 13\}, \{14, 15, 16\}\}.
  \end{eqnarray*}
  The 4rth aggregate $\{10, 11\}$ are denoted by solid spheres in the
  left sub figure of Figure (\ref{aggregation}).
\item For MSCN method, we identify the nodes in graph of $\mathscr{D}$ which have direct
  edges to the nodes of the aggregate $\{10, 11\}$. The identified
  nodes are $\{10, 11\}$, i.e., the nodes shaded with horizontal line
  patterns.
\item Finally assemble the following matrix
$$
 \mathbf{C_4}=\left[\begin{array}{cc|cc}
    D(10,10) & D(10,11)  & F(10,10) & F(10,11)  \\
    D(11,10) & D(11,11) & F(11,10) & F(11,11) \\ \hline
    E(10,10) & E(10,11)  & G(10,10) & G(10,11)  \\
    E(11,10) & E(11,11)  & G(11,10) & G(11,11)  \\
  \end{array}
\right]$$ \normalsize The matrix above is presented with entries with
local numbering. The assembled matrix in terms of the global numbering
is the following
$$
 \mathbf{C_4}=\left[\begin{array}{cc|cc}
    C(10,10) & C(10,11)  & C(10,26) & C(10,27)  \\
    C(11,10) & C(11,11) & C(11,26) & C(11,27) \\ \hline
    C(26,10) & C(26,11)  & C(26,26) & C(26,27)  \\
    C(27,10) & C(27,11)  & C(27,26) & C(27,27)  \\
  \end{array}
\right].$$ \normalsize

\item Denote the above block $2 \times 2$ matrix as follows
$$
\mathbf{C_{4}}=\left[\begin{array}{c|c} D_{4} & E_{4} \\ \hline F_{4}
    & G_{4}
  \end{array}
\right].$$ Thus, the mini Schur complement for the $4rth$ aggregate is
given as follows
\begin{eqnarray}
  S_{4} = G_{4} - F_{4}(D_{4})^{-1}E_{4}
\end{eqnarray}
Finally, the columns corresponding to the $4rth$ aggregate, i.e., for the
aggregate $\{10,11\}$ in the mini Schur complement $S_4$ are $S_{4}(:,1)$ and
$S_{4}(:,2)$ in the familiar Matlab notation.  We substitute these
columns in the global Schur complement $S$. In the global numbering,
the columns $S_{4}(1,:)$ and $S_{4}(2,:)$ maps to the column numbers
$S(10,:)$ and $S(11,:)$ of the global Schur complement to be
estimated. Thus we do the following update to the matrix $S$
\begin{eqnarray*}
  S(10,10:11) = S_{4}(1,:)^{T}, \\
  S(11,10:11) = S_{4}(2,:)^{T}. 
\end{eqnarray*}
Note that we have used the convenient Matlab notation above.
\end{enumerate}

\subsection{Matrix view of the MSCN method}

In this section, we shall present a matrix view of the MSCN method.
Although, the graph view above and the matrix view are both same
method.  The matrix view will lead to a clearer picture and
facilitates proving results.

For MSCN method, no reordering or partitioning is required for the
method to work. However, using efficient partitioning techniques
usually leads us to a considerable saving in the flop count and
enhance parallelism. Thus, in what follows, we shall assume that a
suitable partitioning and reordering has been applied using suitable
graph partitioner. From now onwards, we use the same notations for the
sub-matrices $D$ and $E$ of $\mathbf{C}$ in the reordered matrix
$\mathbf{P^{T}CP}$.

As above, we start with the block $2 \times 2$ partitioned system as
follows

$$
\mathbf{C}=\left[\begin{array}{c|c} D & E \\ \hline F & G
  \end{array}
\right]$$

Then we partition the matrix further as follows

$$
\mathbf{C}=\left[\begin{array}{c|c}
    
    \begin{array}{c|c} D_{11} & D_{12} \\ \hline D_{21}
      & D_{22}
    \end{array}

    &       \begin{array}{c|c} E_{11} & E_{12} \\ \hline E_{21}
      & E_{22}
    \end{array} \\ \hline

    \begin{array}{c|c} F_{11} & F_{12} \\ \hline F_{21}
      & F_{22}
    \end{array}     &         \begin{array}{c|c} G_{11} & G_{12} \\ \hline G_{21}
      & G_{22}
    \end{array}

  \end{array}
\right]$$

Now we construct a sparse approximation of matrix $\mathbf{C}$ above
by dropping the blocks $D_{ij},E_{ij},F_{ij},$ and $G_{ij}$ for which
$i \neq j$.  As a result, we obtain a first level sparse approximation
denoted by $\mathbf{\widehat{C}_{2}}$ as follows

$$
\mathbf{\widehat{C}_{2}}=\left[\begin{array}{c|c}
    
    \begin{array}{c|c} D_{11} &  \\ \hline 
      & D_{22}
    \end{array}

    &       \begin{array}{c|c} E_{11} &  \\ \hline 
      & E_{22}
    \end{array} \\ \hline

    \begin{array}{c|c} F_{11} &  \\ \hline 
      & F_{22}
    \end{array}     &         \begin{array}{c|c} G_{11} &  \\ \hline 
      & G_{22}
    \end{array}

  \end{array}
\right].$$

Here the subscript 2 in $\mathbf{\widehat{C}_{2}}$ denotes the number
of principle sub matrices of the matrix G. In this case, the two
principle sub matrices we have retained are $G_{11}$ and $G_{22}$.

For the second level of partition, we further partition the blocks
$D_{ii}, E_{ii}, F_{ii}$, and $G_{ii}$ to get a sparse matrix of
$\mathbf{C}$ as follows

$$
\left[\begin{array}{c|c}
    
    \begin{array}{c|c} 

      \begin{array}{c|c} 

        D_{11} & 
      
        D_{12} \\ \hline 
      
        D_{21} & 
  
        D_{22}
      
    \end{array}     
    & 
      

      
      
  
      
    \\

    \hline
      

      
      
  
      
    &
  
    \begin{array}{c|c}

      D_{33} & 
      
      D_{34} \\ \hline 
      
      D_{43} & 
  
      D_{44}
      
    \end{array}

  \end{array}

  &
  \begin{array}{c|c}
      
    \begin{array}{c|c} 

      E_{11} & 
      
      E_{12} \\ \hline 
      
      E_{21} & 
  
      E_{22}
      
    \end{array}     
    & 


      
      
  
      
    \\ \hline


      
      
  
      
    &

    \begin{array}{c|c}

      E_{33} & 
      
      E_{34} \\ \hline 
      
      E_{43} & 
  
      E_{44}
      
    \end{array}     

  \end{array} \\ \hline

  \begin{array}{c|c}

    \begin{array}{c|c} 

      F_{11} & 
      
      F_{12} \\ \hline 
      
      F_{21} & 
  
      F_{22}
      
    \end{array}     
    & 


      
      
  
      
    \\ \hline


      
      
  
      
    &

    \begin{array}{c|c}

      F_{33} & 
      
      F_{34} \\ \hline 
      
      F_{43} & 
  
      F_{44}
      
    \end{array}     

    \end{array}     
    &
    \begin{array}{c|c}

      \begin{array}{c|c} 

        G_{11} & 
      
        G_{12} \\ \hline 
      
        G_{21} & 
  
        G_{22}
      
    \end{array}     
    & 


      
      
  
      
    \\ \hline


      
      
  
      
    &

    \begin{array}{c|c}

      G_{33} & 
      
      G_{34} \\ \hline 
      
      G_{43} & 
  
      G_{44}
      
    \end{array}     

  \end{array}

  \end{array}
\right].$$

Now as before, we construct a sparse approximation of the above matrix
by dropping the blocks $D_{ij},E_{ij},F_{ij}$, and $G_{ij}$ for which
$i \neq j$.  We obtain a second level sparse approximation denoted by
$\mathbf{\widehat{C}_{4}}$ as follows

$$
\mathbf{\widehat{C}_{4}}=\left[\begin{array}{c|c}
    
    \begin{array}{c|c} 

      \begin{array}{c|c} 

        D_{11} & 
      
        \\ \hline 
      
        & 
  
        D_{22}
      
    \end{array}     
    & 
      

      
      
  
      
    \\

    \hline
      

      
      
  
      
    &
  
    \begin{array}{c|c}

      D_{33} & 
      
      \\ \hline 
      
      & 
  
      D_{44}
      
    \end{array}

  \end{array}

  &
  \begin{array}{c|c}
      
    \begin{array}{c|c} 

      E_{11} & 
      
      \\ \hline 
      
      & 
  
      E_{22}
      
    \end{array}     
    & 


      
      
  
      
    \\ \hline


      
      
  
      
    &

    \begin{array}{c|c}

      E_{33} & 
      
      \\ \hline 
      
      & 
  
      E_{44}
      
    \end{array}     

  \end{array} \\ \hline

  \begin{array}{c|c}

    \begin{array}{c|c} 

      F_{11} & 
      
      \\ \hline 
      
      & 
  
      F_{22}
      
    \end{array}     
    & 


      
      
  
      
    \\ \hline


      
      
  
      
    &

    \begin{array}{c|c}

      F_{33} & 
      
      \\ \hline 
      
      & 
  
      F_{44}
      
    \end{array}     

    \end{array}     
    &
    \begin{array}{c|c}

      \begin{array}{c|c} 

        G_{11} & 
      
        \\ \hline 
      
        & 
  
        G_{22}
      
    \end{array}     
    & 


      
      
  
      
    \\ \hline


      
      
  
      
    &

    \begin{array}{c|c}

      G_{33} & 
      
      \\ \hline 
      
      & 
  
      G_{44}
      
    \end{array}     

  \end{array}

  \end{array}
\right]$$

Eliminating the blocks $F_{ii}$ using $D_{ii}$ as pivots, we obtain an
approximation to the original Schur complement $S$ by
$\widehat{S}_{4}=blkDiag(S_{ii})$ as follows
\begin{eqnarray}
  S_{ii} = G_{ii} - F_{ii}D_{ii}^{-1}E_{ii},~i = 1:4
\end{eqnarray}
where the subscript $4$ in $\widehat{S}_{4}$ is the number of
principle sub matrices of matrix G and we call $S_{ii}$ a MSC.
\begin{remark}
  For simplicity, our approach was to partition the matrix recursively
  into block $2 \times 2$ matrix. We could directly identify the
  blocks $G_{ii}$ such that the following expression
  \begin{eqnarray} \label{rec_msc} S_{ii} = G_{ii} -
    F_{ii}D_{ii}^{-1}E_{ii}, i = 1:m
  \end{eqnarray}
  makes sense. Here $m$ denote the number of MSCs desired.
\end{remark}
For notational convenience, we denote diagonal blocks of $D$ by
$\widehat{D}_{m}=blkDiag(D_{ii})$, similarly, we denote
$\widehat{E}_{m}=blkDiag(E_{ii})$, $\widehat{F}_{m}=blkDiag(F_{ii})$,
and $\widehat{G}_{m}=blkDiag(G_{ii})$. Thus, the matrix
$\mathbf{\widehat{C}_{m}}$ in the general case is given as follows
$$
\mathbf{\widehat{C}_{m}} =\left[\begin{array}{c|c} \widehat{D}_{m} &
    \widehat{E}_{m} \\ \hline \widehat{F}_{m} & \widehat{G}_{m}
  \end{array}
\right]$$ Thus, we have $\widehat{S}_{m} = \widehat{G}_{m} -
\widehat{F}_{m}(\widehat{D}_{m})^{-1}\widehat{E}_{m}$ which is an
approximation to the original Schur complement $S=G-FD^{-1}E$.  When
it is not necessary, we shall omit the subscript $m$ from
$\widehat{S_{m}}$.

We are now in a position to formally define the MSCN preconditioner.
\begin{definition}
  Given a block $2\times2$ partitioned matrix $\mathbf{C}$ as follows
$$
\mathbf{C}=\left[\begin{array}{c|c} D & E \\ \hline F & G
  \end{array}
\right]$$ The MSCN preconditioner denoted by $\mathbf{B_{MSCN}}$ is
defined as follows
\begin{equation}
  \mathbf{B_{MSCN}}= 
  \left[
    \begin{array}{cc}
      D &   \\
      F & \widehat{S} \\
    \end{array} 
  \right]
  \left[
    \begin{array}{cc}
      D^{-1} &      \\
      & (\widehat{S})^{-1}  \\
    \end{array}
  \right]  
  \left[   
    \begin{array}{cc}
      D & E   \\
      &  \widehat{S}   \\ 
    \end{array}
  \right],
  \label{NDfactor}
\end{equation} \label{msc_def}
If $m$ is the number of MSCs considered, then $\widehat{S}=blkDiag(S_{11},\dots,S_{mm})$, where $S_{ii}$ is given by the formula (\ref{rec_msc}) above.
\end{definition}
\begin{definition}
  The LUM method is defined as above except that
  $\widehat{S}=blkDiag(G_{11},\dots,G_{mm})$. Here $m$ is the number
  of MSCs.
\end{definition}
\begin{theorem}
  If D defined above is symmetric positive definite, then
  $\mathbf{\widehat{S}_{m}}$ exists. Thus, the MSCN preconditioner
  exists.
\end{theorem}
\textbf{Proof:} If the (1,1) block D of the matrix $\mathbf{C}$ is SPD
then each $D_{ii}$ (being the diagonal blocks of D) for all $1 \leq i
\leq m$ are SPD and the formula (\ref{rec_msc}) does not break down.

\begin{lemma} \label{ban} \cite{ban} Given a block $2 \times 2$
  partitioned matrix as follows
$$
\mathbf{Z}=\left[
  \begin{array}{cc}
    D & E                         \\
    F              & G \\
  \end{array}
\right],
$$
where the sub matrix D is nonsingular, then the inverse of the matrix
$\mathbf{Z}$ is given as follows
$$
\mathbf{Z}^{-1}=\left[
  \begin{array}{cc}
    D^{-1}+D^{-1}ES^{-1}FD^{-1} & -D^{-1}ES^{-1}                         \\
    -S^{-1}FD^{-1}              & S^{-1} \\
  \end{array}
\right],
$$
where $S = G-FD^{-1}E$.
\end{lemma}

\textbf{Proof:} The proof seems to appear first in \cite{ban}. The
result is known as Banachiewicz inversion formula for the inverse of a
block $2 \times 2$ partitioned matrix.

\begin{theorem} \label{leftPre} Consider the block $2 \times 2$
  partitioned matrix as follows
$$
\mathbf{C}=\left[
  \begin{array}{cc}
    D_{p \times p} & E                         \\
    F              & G \\
  \end{array}
\right]_{n \times n},
$$
and $\mathbf{B_{MSCN}}$ be the MSCN preconditioner as defined in
Def. (\ref{msc_def}).  Then the MSCN preconditioned matrix
$\mathbf{C}$ i.e., $(\mathbf{B_{MSCN}})^{-1}\mathbf{C}$ has $p$
eigenvalues exactly equal one.  The rest of the n-p eigenvalues are
the eigenvalues of $(\widehat{S}_{m})^{-1}S$, where $S = G - FD^{-1}E$
is the complete Schur complement of the matrix $\mathbf{C}$.
\end{theorem}

\textbf{Proof:} Using Lemma (\ref{ban}) above, we have

\begin{equation}
  \mathbf{B^{-1}_{MSCN}}= 
  \left[
    \begin{array}{cc}
      I_{p \times p}                                      & -D^{-1}E {\widehat{S}_{m}}^{-1} \\
      0                                      & {\widehat{S}_{m}}^{-1}          \\
    \end{array} 
  \right]_{n \times n}
  \left[
    \begin{array}{cc}
      D_{p \times p}                                      &                                           \\
      & \widehat{S}_{m}               \\
    \end{array}
  \right]_{n \times n}  
  \left[   
    \begin{array}{cc}
      D_{p \times p}^{-1}                                 & 0                                         \\
      {\widehat{S}_{m}}^{-1}FD^{-1} & {\widehat{S}_{m}}^{-1}          \\ 
    \end{array}
  \right]_{n \times n},
\end{equation} 

and the MSCN preconditioned matrix is given as follows

\begin{equation} \label{expLeftPre} \mathbf{B_{MSCN}^{-1}C} = \left[
    \begin{array}{cc}
      I_{p \times p}                                 &  D^{-1}E(I - {\widehat{S}_{m}}^{-1}S) \\
      0 & {\widehat{S}_{m}}^{-1}S          \\ 
    \end{array}
  \right]_{n \times n},
\end{equation}
where $S = G - FD^{-1}E$. Hence the theorem.

\begin{remark}
  In theorem (\ref{leftPre}) above, we have proved the result for left
  preconditioned matrix. However, a similar result holds for Right
  preconditioned matrix as follows

  \begin{equation} \label{expRightPre} \mathbf{CB_{MSCN}^{-1}} =
    \left[
      \begin{array}{cc}
        I_{p \times p}                                 & 0 \\
        (I - S{\widehat{S}_{m}}^{-1})FD^{-1} & S{\widehat{S}_{m}}^{-1}          \\ 
      \end{array}
    \right]_{n \times n},
  \end{equation}
  where $S = G - FD^{-1}E$.
\end{remark}

\begin{remark}
  Comparing the expression for the left preconditioned matrix in
  Equation (\ref{expLeftPre}) with that of right preconditioned matrix
  in Equation (\ref{expRightPre}), we find that
  $S{\widehat{S}_{m}}^{-1}$ and ${\widehat{S}_{m}}^{-1}S$ are similar
  matrices. Thus, we have that the n-p eigenvalues for left
  preconditioned matrix given by eigenvalues of
  $S{\widehat{S}_{m}}^{-1}$ are same as the n-p eigenvalues of the
  right preconditioned matrix given by eigenvalues of
  ${\widehat{S}_{m}}^{-1}S$. In practice, the original matrix
  $\mathbf{C}$ is indefinite as confirmed in top left figure in Table
  (\ref{spec}) where we notice many negative eigenvalues.
  \begin{table}
    \caption{\label{spec} Top Left: Real part of the eigenvalues of original matrix, Top Right
      Real part of the eigenvalues of the MSCN right preconditioned matrix, Bottom left: Eigenvalues (Real and imaginary) of original matrix, Bottom right: Eigenvalues (Real and imaginary) of MSCN right preconditioned matrix}
    \begin{tabular}{ll}
      \includegraphics[scale=0.35]{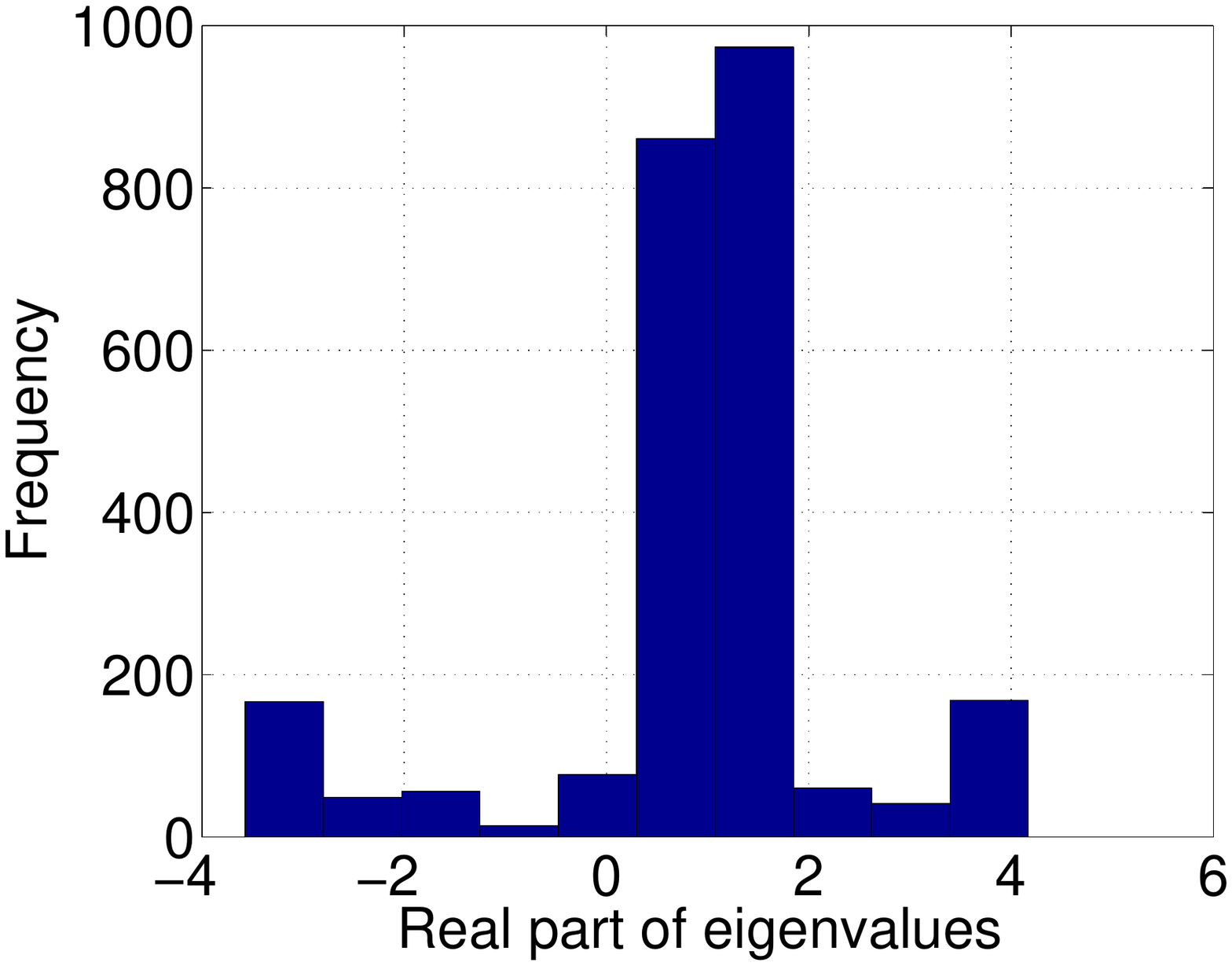} 
      \label{histplotC}
      &
      \includegraphics[scale=0.35]{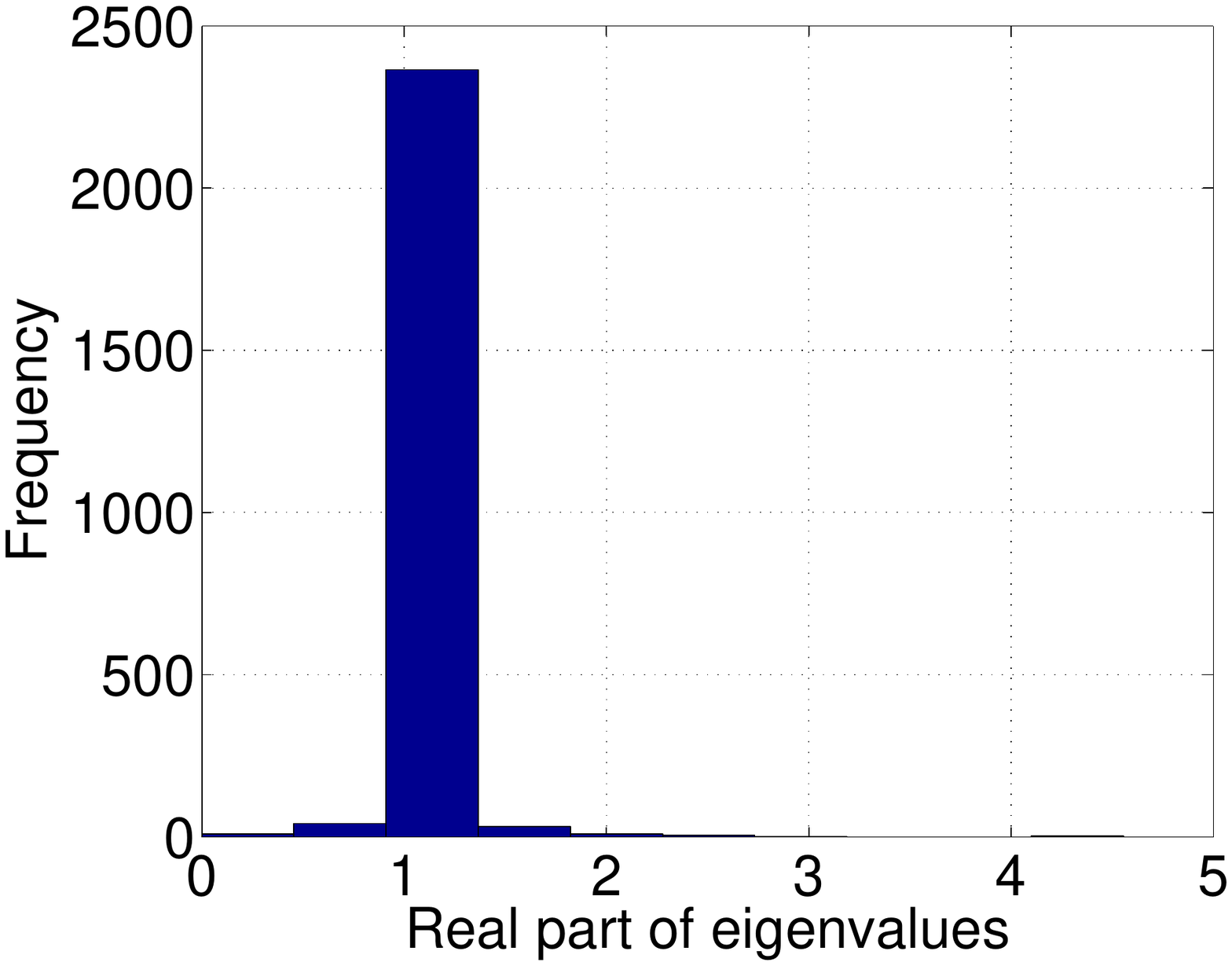}  \\

      \includegraphics[scale=0.35]{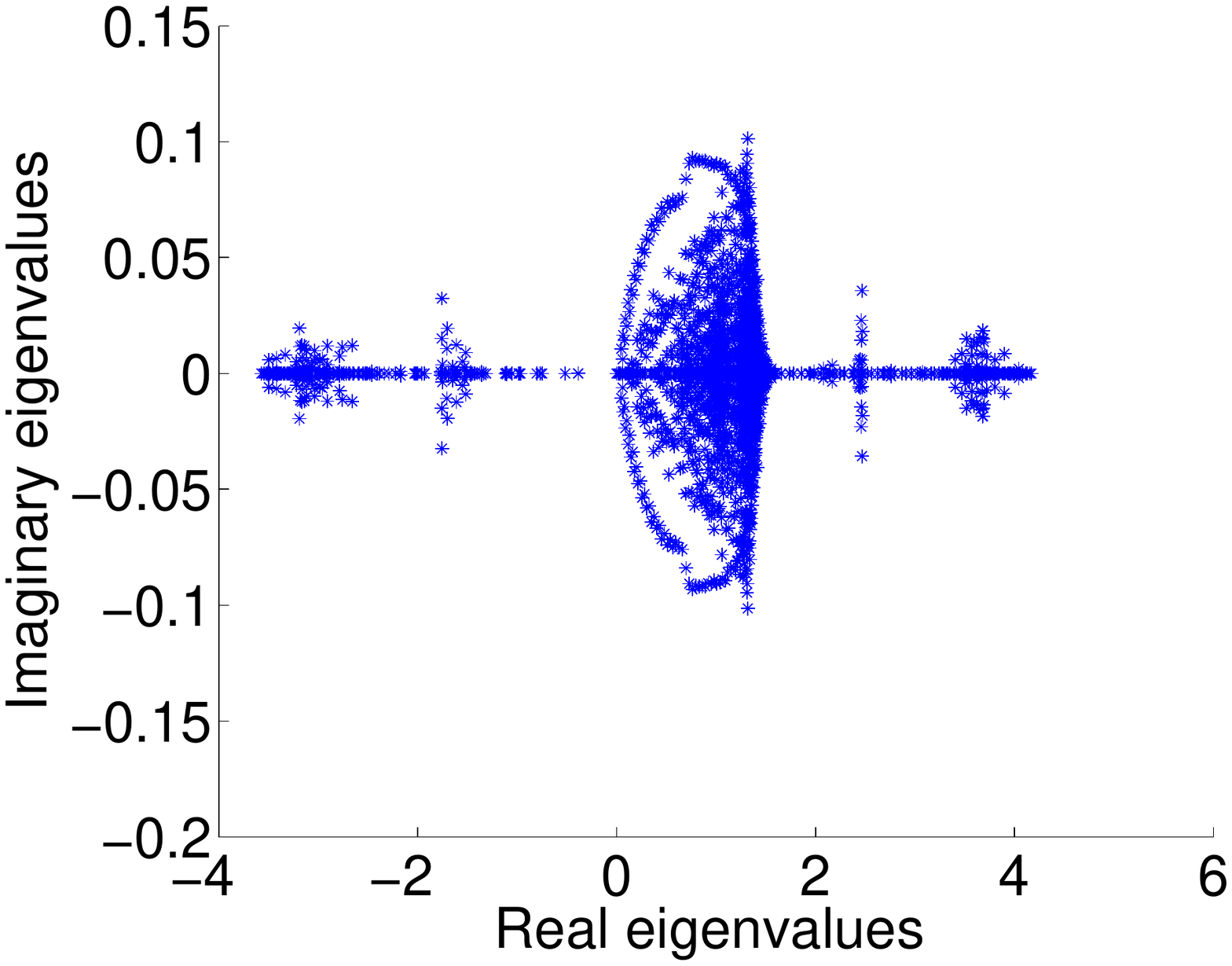} 
      \label{histplotC}
      &
      \includegraphics[scale=0.35]{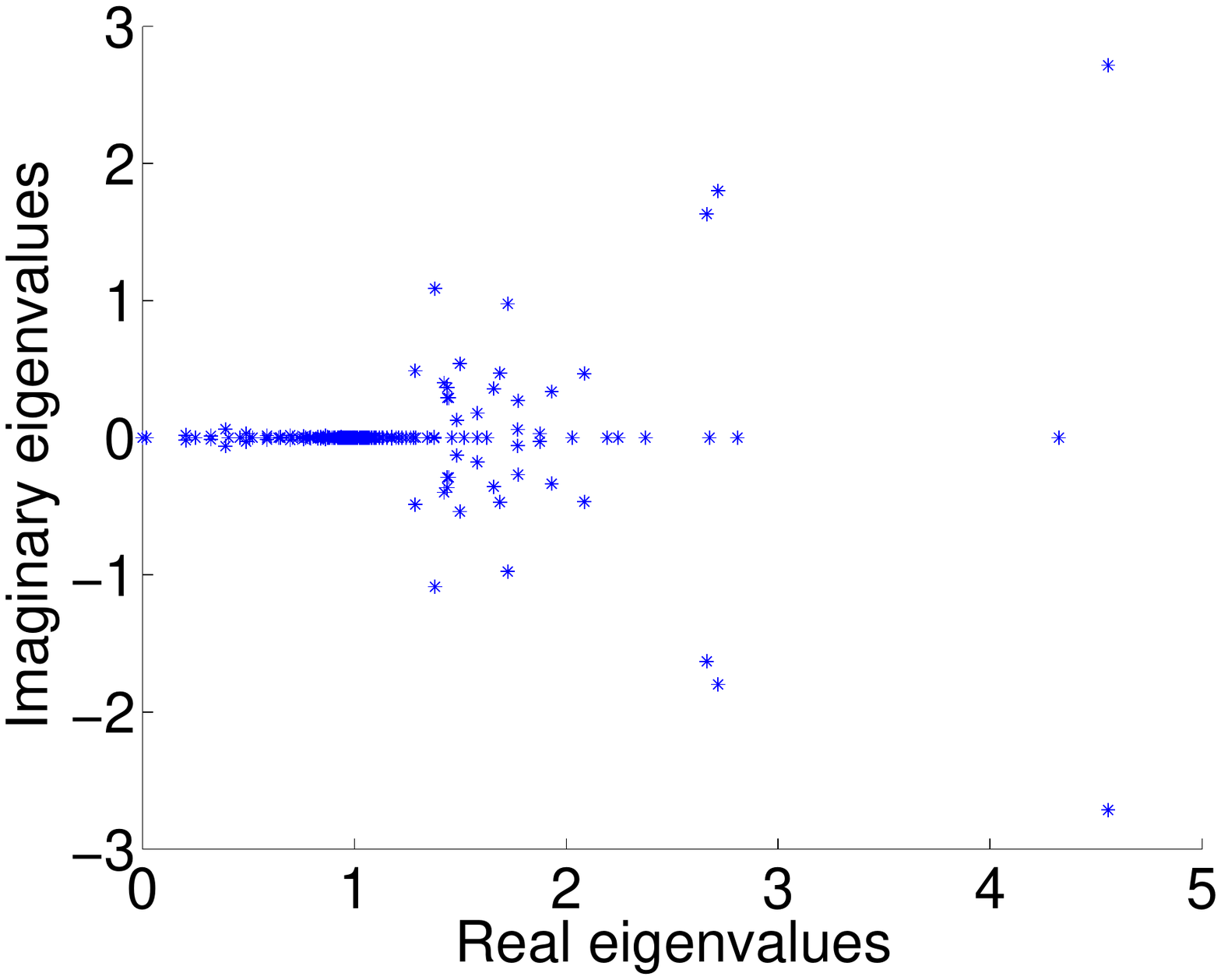}  \\
  
    \end{tabular}
  \end{table}

\end{remark}

\begin{remark}
  In general, it is difficult to estimate the eigenvalues of
  ${\widehat{S}_{m}}^{-1}S$ since we need some additional assumptions
  for the sub matrices of $E$ and $F$ which, in practice, are
  unknown. However, it may be a good idea to partition the matrix
  which leads to a large (1,1) block D, and as a consequence more and
  more eigenvalues are equal to one. On the other hand, since we will
  need to solve equations of the form $Dx=y$, so ideally D should be
  partitioned into small diagonal blocks that are easier to invert.
  Thanks to several graph partitioning softwares readily available, we
  could easily obtain such reordering.  Some of the popular
  partitioning and reordering softwares available in the public domain
  are namely, METIS \cite{kar3} and independent set ordering used in
  \cite{saad2002, saad2003}. Use of such partitioning and reordering
  techniques leads to a purely algebraic domain decomposition which
  takes matrix as an input and partitions the graph (rather than the
  computational domain) by selecting a separator (set of edges or
  vertices), removal of which leads to several disconnected subgraphs.
\end{remark}

The quality of the preconditioner is determined by the quality of the
Schur complement approximation and the spectrum distribution of
$\widehat{S}^{-1}S$.

Consider the block $2 \times 2$ partitioned matrix as follows
$$
\mathbf{C}=\left[
  \begin{array}{cc}
    D & E                         \\
    F              & G \\
  \end{array}
\right]_{n \times n},
$$

where $D = blkDiag(D_{1},D_{2})$,

\begin{table}[h]
  \begin{tabular}{lll}
    $  E= \left[
      \begin{array}{cc}
        E_{11} & E_{12}                         \\
        E_{21} & E_{22} \\
      \end{array}
    \right]_{n \times n},$
    &
    $  F=\left[
      \begin{array}{cc}
        F_{11} & F_{12}                         \\
        F_{21} & F_{22} \\
      \end{array}
    \right]_{n \times n},$
    &
    $    G=\left[
      \begin{array}{cc}
        G_{11} & G_{12}                         \\
        G_{21} & G_{22} \\
      \end{array}
    \right]_{n \times n},$
  \end{tabular}
\end{table}

The Schur complement $S$ is given as follows

\begin{equation}
  S=\left[
    \begin{array}{cc}
      G_{11}-E_{11}D_{1}^{-1}F_{11}-E_{12}D_{2}^{-1}F_{21} & G_{12}-E_{11}D_{1}^{-1}F_{12}-E_{12}D_{2}^{-1}F_{22}                         \\
      G_{21}-E_{21}D_{1}^{-1}F_{11}-E_{22}D_{2}^{-1}F_{21} & G_{22}-E_{21}D_{1}^{-1}F_{12}-E_{22}D_{2}^{-1}F_{22} \\
    \end{array}
  \right]
\end{equation} 

The approximation to the Schur complement via MSCN method for this
model problem is given as follows
\begin{equation}
  \widehat{S}=\left[
    \begin{array}{cc}
      G_{11}-E_{11}D_{1}^{-1}F_{11} &                               \\
      & G_{22}-E_{21}D_{1}^{-1}F_{12} \\
    \end{array}
  \right].
\end{equation}   
Thus, $\widehat{S}$ is an approximation to the block Jacobi
preconditioner for $S$. Thus, for the general case with many
subdomains, $\widehat{S}$ remains a crude approximation to
the block Jacobi preconditioner of $S$.

\subsection{Previous work}
The methods presented in this work are related to patch method as
introduced in \cite{rou1,rou2}. However, there are significant
differences. We list these differences as follows

\begin{enumerate}
\item First, we start with an aggregate, i.e., a group of nodes, whereas,
  in \cite{rou1,rou2} the Schur complement is built around a node. For
  a sufficiently large aggregation, the columns corresponding to the
  whole group of nodes of an aggregate are approximated, resulting in a
  much faster sweep in the graph.
\item The aggregation scheme for the MSCN method considered in this
  paper is based on the numbering of the nodes in the grid. On the
  other hand, in the patch method, a small patch is formed around a
  node based on the edge connections. For example, in Figure
  (\ref{aggregation}), the nearest neighbors of node numbered
  12 are nodes numbered 11 and 13. But for the MSCN method, the nearest
  neighbors of 12 would be nodes numbered 13 and 14.
\end{enumerate}

\section{Exploiting parallelism and implementation aspects}
The methods we have seen in the previous section are massively
parallel both during the setup phase as well as during the iteration
phase. In this section, we understand the parallelism in the
method. In Algorithm (\ref{pseudocode_build}), we provide the
pseudocode for the construction phase of the MSC methods and in
Algorithm (\ref{pseudocode_solve}) the solution procedure is
presented.

\begin{algorithm} 
  \caption{\label{pseudocode_build} PSEUDOCODE TO CONSTRUCT $B$ FOR
    MSCN,MSCNR,LUM,OMSCN,OMSCNR, and OLUM methods} 
\begin{algorithmic}
\STATE INPUT: 
\begin{itemize}
 \item $$
   \mathbf{C}=\left[\begin{array}{c|c} D & E \\ \hline F & G
     \end{array}
   \right]$$

 \item $k$ = Number of MSCs desired 
 \item $w$ = array of length $k$ contains size overlap, $w_{k}=0$
\end{itemize}

\STATE // Find principle submatrices (aggregates)
\IF{MSCN, LUM, or MSCNR}
  \STATE Find principle submatrices of $G$, $D$, $F$, and $E$ as follows
  \STATE $P_{G}$ = $\{G_{11},G_{22},\dots,G_{kk}\}$, $G_{ii}=G(r_{i}:r_{i+1}-1;r_{i}:r_{i+1}-1)$
  \STATE $P_{D}$ = $\{D_{11},D_{22},\dots,D_{kk}\}$, $D_{ii}=D(r_{i}:r_{i+1}-1;r_{i}:r_{i+1}-1)$
  \STATE $P_{F}$ = $\{F_{11},F_{22},\dots,F_{kk}\}$, $F_{ii}=F(r_{i}:r_{i+1}-1;r_{i}:r_{i+1}-1)$
  \STATE $P_{E}$ = $\{E_{11},E_{22},\dots,E_{kk}\}$, $E_{ii}=E(r_{i}:r_{i+1}-1;r_{i}:r_{i+1}-1)$
  \STATE Here, $r_{1}=1$, $r_{k+1}= ncols(G)+1$
\ELSIF{OMSCN, OLUM, OMSCNR}
  \STATE Find principle submatrices of $G$, $D$, $F$, and $E$ as follows
  \STATE $P_{G}$ = $\{G_{11},G_{22},\dots,G_{kk}\}$, $G_{ii}=G(r_{i}:r_{i+1}-1+w_{i};r_{i}:r_{i+1}-1+w_{i})$
  \STATE $P_{D}$ = $\{D_{11},D_{22},\dots,D_{kk}\}$, $D_{ii}=D(r_{i}:r_{i+1}-1+w_{i};r_{i}:r_{i+1}-1+w_{i})$
  \STATE $P_{F}$ = $\{F_{11},F_{22},\dots,F_{kk}\}$, $F_{ii}=F(r_{i}:r_{i+1}-1+w_{i};r_{i}:r_{i+1}-1+w_{i})$
  \STATE $P_{E}$ = $\{E_{11},E_{22},\dots,E_{kk}\}$, $E_{ii}=E(r_{i}:r_{i+1}-1+w_{i};r_{i}:r_{i+1}-1+w_{i})$
  \STATE Here, $r_{1}=1$, $r_{k+1}= ncols(G)+1$, $w_{i}< (r_{i+1}-r_{i})$
\ENDIF
 \STATE // Now construct MSCs
\IF{MSCN or OMSCN}

\FOR {$i$ to $k$}
\STATE $S_{ii} = G_{ii} - F_{ii}D_{ii}^{-1}E_{ii}$ 
\ENDFOR

\ELSIF{LUM}

\FOR {$i$ to $k$}
\STATE $S_{ii} = G_{ii}$  
\ENDFOR

\ELSIF{OMSCNR or MSCNR}

\FOR {$i$ to $k$}
\STATE $S_{ii} = G_{ii} - F_{ii}D_{ii}^{-1}(E_{ii}\mathbf{1})$ 
    \STATE where $\mathbf{1}=[1,1,\dots,1]^{T}$
\ENDFOR
\ENDIF

\STATE Set $\widehat{S}=blkDiag(S_{11},S_{22},\dots,S_{kk})$

\STATE OUTPUT: 
\begin{equation}
  \mathbf{B}= 
  \left[
    \begin{array}{cc}
      D &   \\
      F & \widehat{S} \\
    \end{array} 
  \right]
  \left[ 
    \begin{array}{cc}
      D^{-1} &      \\
      & (\widehat{S})^{-1}  \\
    \end{array}
  \right]  
  \left[   
    \begin{array}{cc}
      D & E   \\
      &  \widehat{S}   \\ 
    \end{array}
  \right],
\end{equation}
\STATE Here B is either of the MSC based preconditioners
\end{algorithmic}
\end{algorithm}

\begin{algorithm} 
  \caption{\label{pseudocode_solve} PSEUDOCODE TO SOLVE WITH $B$ FOR
    MSCN,MSCNR,LUM,OMSCN,OMSCNR, and OLUM methods} 
\begin{algorithmic}

\STATE OBJECTIVE
\begin{itemize}
\item To solve
\begin{equation}
  \mathbf{B}= 
  \left[
    \begin{array}{cc}
      D &   \\
      F & \widehat{S} \\
    \end{array} 
  \right]
  \left[ 
    \begin{array}{cc}
      D^{-1} &      \\
      & (\widehat{S})^{-1}  \\
    \end{array}
  \right]  
  \left[   
    \begin{array}{cc}
      D & E   \\
      &  \widehat{S}   \\ 
    \end{array}
  \right]
  \left[   
    \begin{array}{c}
        x_{1}    \\
        x_{2}    \\ 
    \end{array}
  \right]=
  \left[   
    \begin{array}{c}
        y_{1}    \\
        y_{2}    \\ 
    \end{array}
  \right]
\end{equation}
\item $D=blkDiag(A_{11},A_{22},\dots,A_{mm})$, $\widehat{S}=blkDiag(S_{11},S_{22},\dots,S_{kk})$
\item $m$ := number of partitions, $k$ := number of MSCs, $y$ := Right hand side  
\end{itemize}
\STATE // forward sweep
\FOR{i= 1 to $m$}
\STATE $z_{1i} = A_{ii}^{-1}(y_{1i})$. Here $z_{1}=\left[z_{11},\dots,z_{1m}\right]$, $y_{1}=\left[y_{11},\dots,y_{1m}\right]$ // Can use an inexact solve
\ENDFOR
\FOR{i=1 to $k$}
\STATE $z_{2i} = S_{ii}^{-1}(y_{2i}-(Ez_{1})_{i})$. Here $z_{2}=\left[z_{21},z_{22},\dots,z_{2k}\right]$, $y_{2}=\left[y_{21},y_{22},\dots,y_{2k}\right]$. // Can use an inexact solve
\ENDFOR
\STATE // backward sweep
\STATE Set $x_{2}=z_{2}$
\FOR{i=1 to $m$}
\STATE $x_{1m} = z_{1m}-A_{ii}^{-1}(Fx_{2})_{i}$, $x_{1} = \left[x_{11},x_{12},\dots,x_{1m}\right]$. // Can use an inexact solve
\ENDFOR
\STATE OUTPUT: $x$
\end{algorithmic}
\end{algorithm}

We start with the aggregation process by selecting a set of aggregated
nodes, and subsequently we construct mini Schur complement for each
aggregate. The construction of each of the mini Schur complements are
independent of each other. Thus, if $N$ is the time required to
compute the approximated global Schur complement sequentially, then
since each of the mini Schur complements could be computed by a
processor, ideally, we have a decrease in the computation time for the
setup phase of the preconditioner by $N/p$, $p$ being the number of processors.

Obviously, the size of the MSCs should be large enough and should
involve significant computational work compared to the overhead
involved in setting up the parallel task.

\begin{figure}
  \caption{Scalability curve for the construction of MSCN for a $32
    \times 32$ grid leaky lid driven cavity problem on four cores
    using parallel computing toolbox of Matlab 7.10}
  \label{scal}
  \begin{center}
    \includegraphics[scale=0.4]{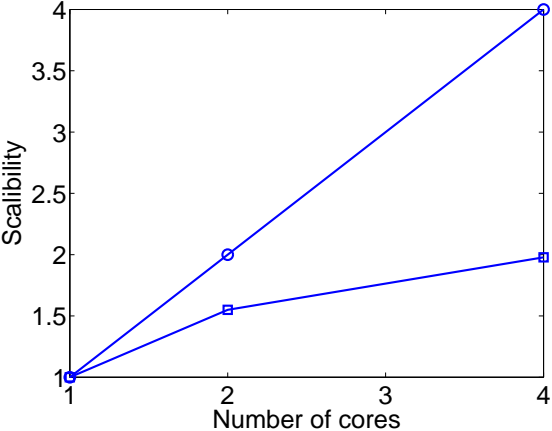} \\
 \end{center}
\end{figure}

For MSCN and MSCE methods, the global Schur complements are block
diagonal matrices, see left figure in Table (\ref{matrix_no_overlap}). Thus, the
factorization phase of the approximated global Schur complement is also completely
parallel. Moreover, due to the same reasons, the solve phase with
approximate Schur complement can de done in parallel. In Figure
(\ref{scal}), we show the scalability in the construction phase of the
MSCN preconditioner for a leaky lid driven cavity problem on $32
\times 32 $ grid. We observe that on four cores, the construction
phase has a speedup of about two. The scalability result is obtained
with parallel computing toolbox of Matlab. The parallel program for
construction is easily implemented by replacing the keyword ``for'' by
the keyword ``parfor'' as follows
\begin{eqnarray*}
  \text{parfor}~ i &=& 1 ~\text{to number of MSCs} \\
  S_{ii} &=& G_{ii} - F_{ii}D_{ii}^{-1}E_{ii} \\
  \text{end}~~~~~&~& \\
\end{eqnarray*}
Here the sub matrix $D_{ii}$ is reasonably small to be inverted easily
by decomposing it into exact triangular factors: $\left[L_{ii},U_{ii}\right] =
lu(D_{ii}).$ While computing $D_{ii}^{-1}E_{ii}$, we achieve one more
level of parallelism by having factorized the matrix $D_{ii}$ into a
product of lower and upper triangular factors, and then solving with
column of $E_{ii}$ as the right hand side. On the other hand, in the
overlapping MSC case, although the computation of MSC can be done
independently, the resulting approximation to the global Schur
complement is not block diagonal as seen in the right figure in Table
(\ref{matrix_no_overlap}). 

  \begin{table}
    \caption{ \label{matrix_no_overlap}Left: structure of approximated
      Schur complement $\mathbf{\widehat{S}_{5}}$ for MSCN, LUM, and
      MSCNC. Right: structure of the approximated Schur complement
      $\mathbf{\widehat{S}_{5}}$ for OMSCN,OMSCNC, and OLUM.}
    \begin{tabular}{ll}
      \includegraphics[scale=0.35]{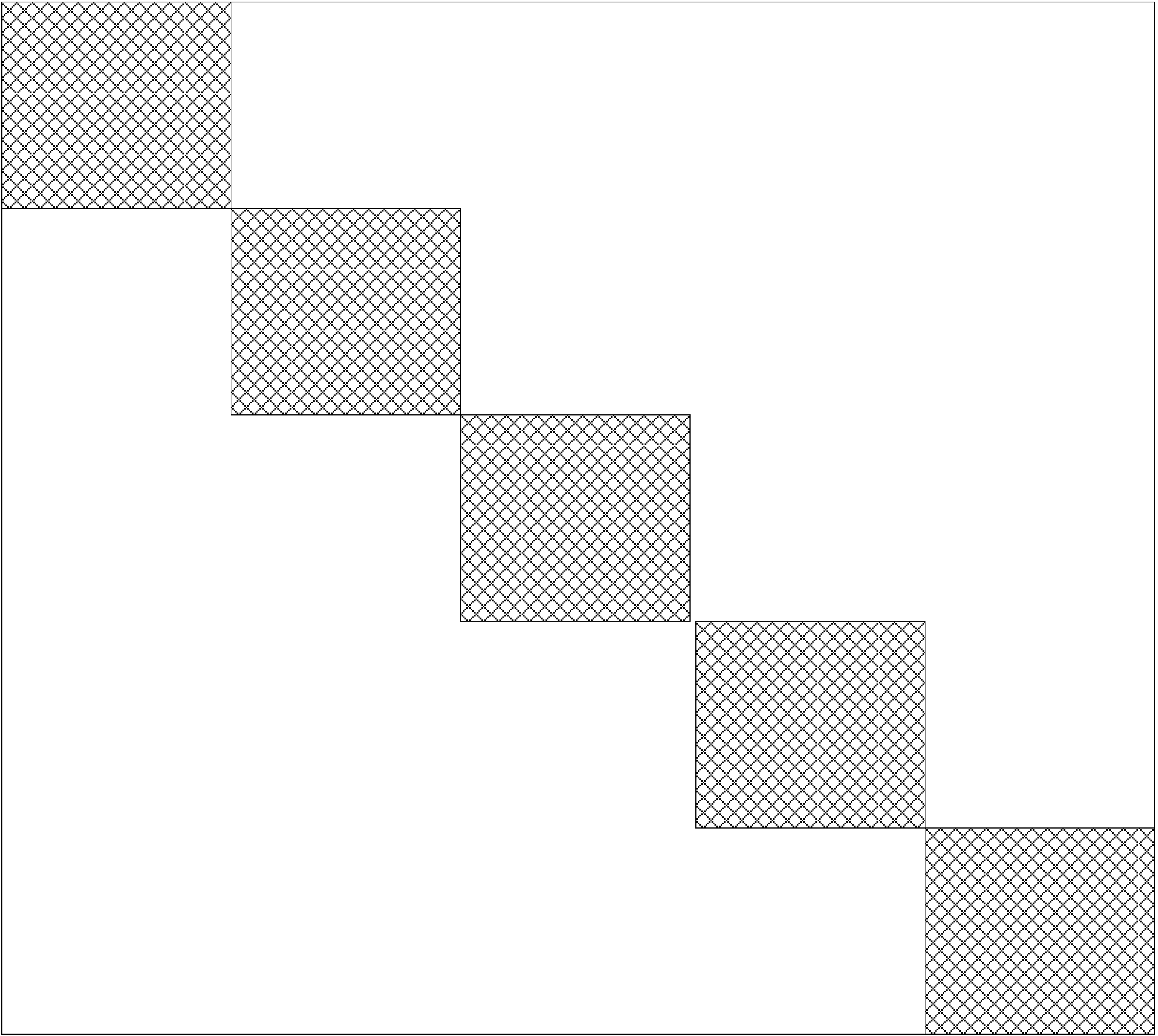}
      &
    \includegraphics[scale=0.35]{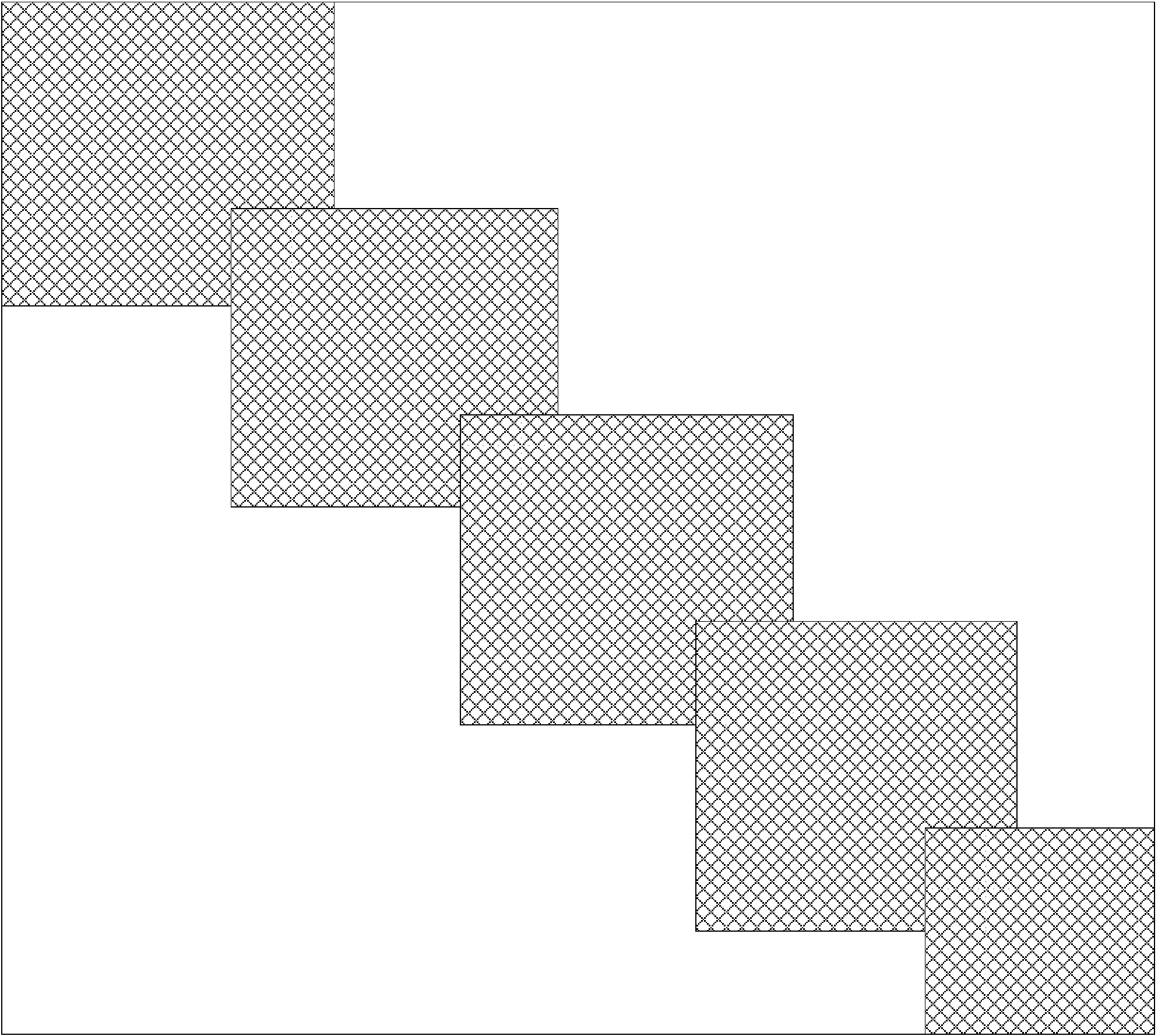} 
    \end{tabular}
  \end{table}


\section{Numerical experiments}
\begin{figure}
  \caption{\label{Grid1} Q2-Q1 discretization for leaky lid driven
    cavity for h=1/128 stretched grid (stretch factor=1.056) with
    Reynolds number of 1000}
  \begin{center}
    \includegraphics[scale=0.3]{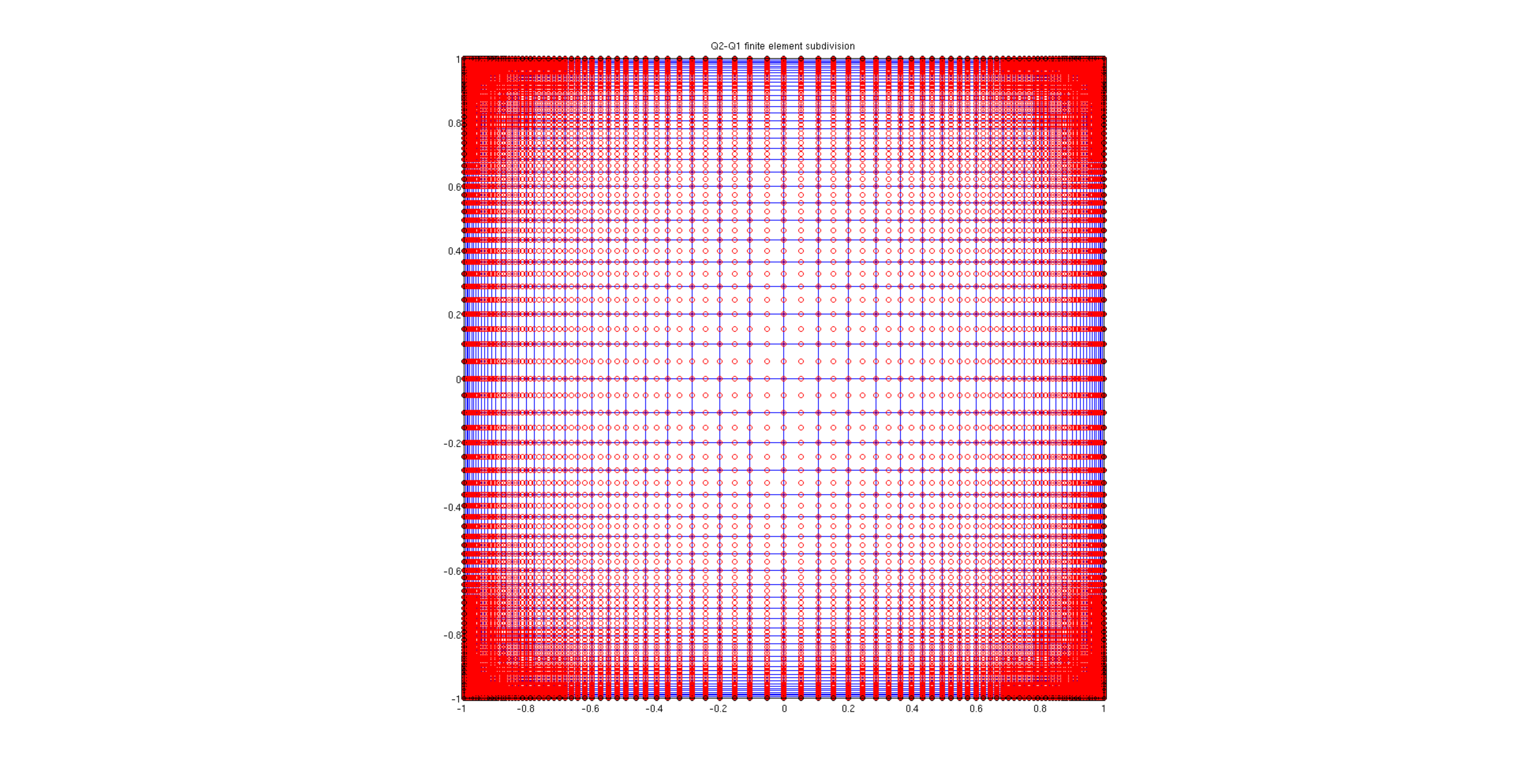} \\
  \end{center}
\end{figure}

The numerical experiments were performed on Matlab 7.10 in double
precision arithmetic with multi-threading enabled on a core i7 (720QM)
Intel processor with 8 processing threads. The system had 6GB of DDR3
RAM and 6MB of L3 cache. The iterative accelerator used is restarted
GMRES with subspace dimension 300. We keep the subspace dimension
large to keep avoid the effects of restart. The maximum number of
iterations allowed is 3000 and the stopping criteria is the decrease
of relative residual below $10^{-9}$. The given coefficient matrix is
scaled by dividing each row by an entry on the same row with maximum
absolute value. For the sake of comparison with sequential PCD and LSC
methods, the experiments with the new methods are also done
sequentially.  The test set consists of standard leaky lid driven
cavity problem defined on both uniform and stretched grid. The
discretization scheme used is the Q2-Q1 (biquadratic velocity/bilinear
pressure) mixed finite element discretization with viscosity varying
from 0.1 to 0.001. A sample of the discretization of Q2-Q1 scheme for
$h=1/64$ for stretched grid with stretch factor of 1.056 with more
finer discretization near the boundaries is shown in figure (\ref{Grid1}).
These problems are generated using the IFISS software \cite{elm}. We
compare all three important aspects of the methods, namely, the
storage requirements, i.e., the fill factor, the iteration count, and
the CPU time.  For MSCE method, the parameters $p_{i}$ and $r_{i}$ are
taken to be equal to ``sz'' defined in Table (\ref{not}).  \small
\begin{table}
  \caption{Notations used in tables of numerical experiments }
  \label{not}
  \begin{center}
    \begin{tabular}{ll}
      \hline
      Notations & Meaning                                                           \\
      \hline
      grid      & Number of discretization points                                   \\
      sz        & Size of the Mini Schur complement                                 \\
      nA        & Number of independent blocks in the (1,1) block                   \\
      nS        & Number of MSCs considered                                         \\
      sA        & Threshold size of the diagonal blocks of (1,1) block              \\
      tolA      & Tolerance for the incomplete LU approx. for (1,1) block           \\
      its       & Iteration count                                                   \\
      time      & time for construction and solution excludes partition time        \\
      MSCN      & Mini Schur complement with numbering based aggregation            \\
      OMSCN     & Overlapped mini Schur complement with numbering based aggregation \\
      LUM       & Lumped approximation                                              \\
      MSCE      & Mini Schur complement with edge based aggregation                 \\
      MPCD      & Modified pressure convection diffusion                            \\
      LSC       & Least square commutator                                           \\
      NA        & Not applicable                                                    \\
      NC        & Not converged                                                     \\
      -         & Test abandonned due to relatively large time                      \\
      \hline
    \end{tabular}
  \end{center}
\end{table}
\normalsize
\subsection{Leaky lid driven cavity}
\begin{figure}
  \caption{\label{stream2} Left: Velocity streamlines, Right: Pressure
    field for leaky lid driven cavity for h=1/128 stretched grid with
    Q2-Q1 discretization with Reynolds number of 1000}
  \begin{center}
    \includegraphics[scale=0.3]{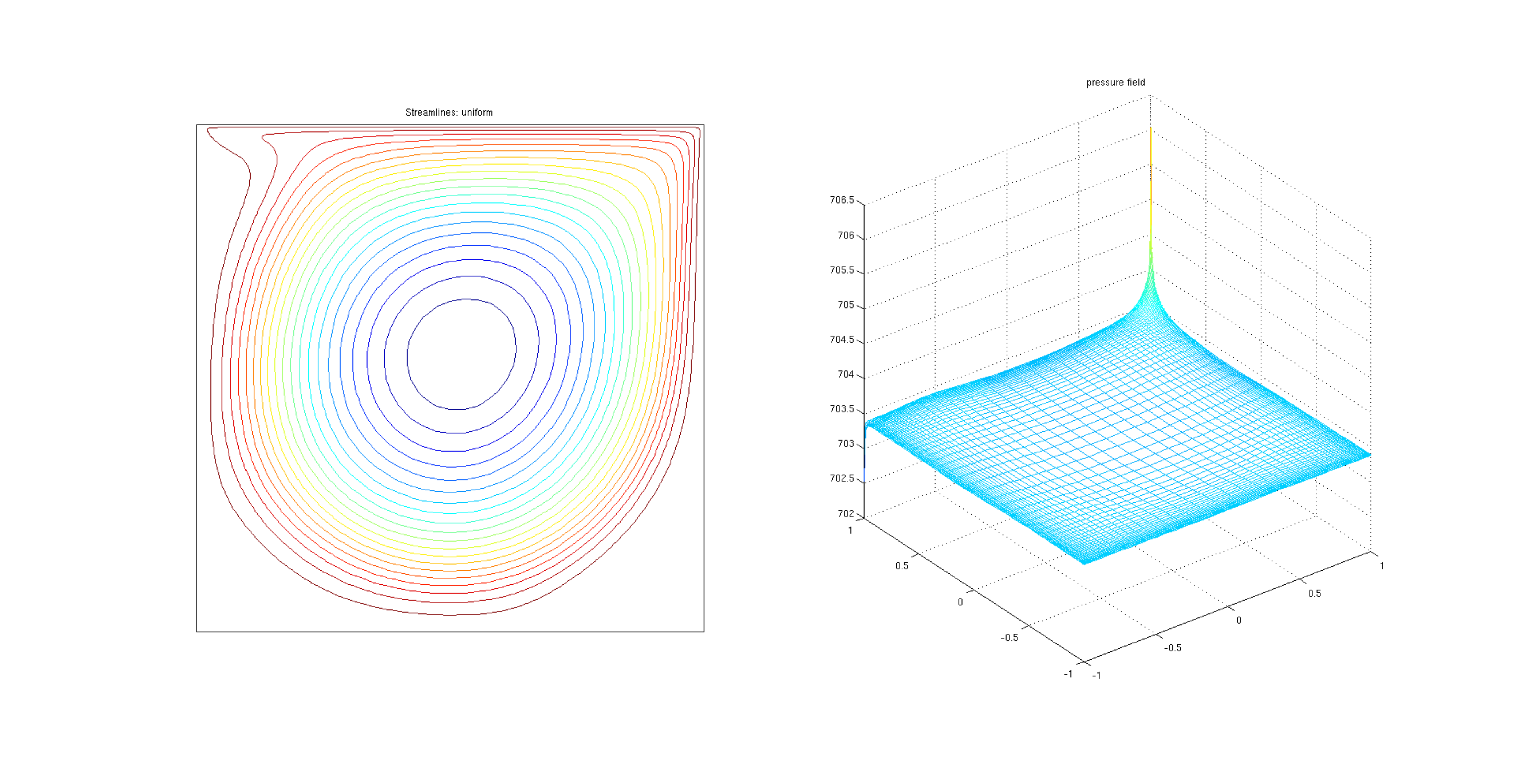} \\
  \end{center}
\end{figure}

\begin{figure}
  \caption{\label{stream1} Left: Velocity streamlines, Right: Pressure
    field for leaky lid driven cavity for h=1/128 stretched grid with
    Q2-Q1 discretization with Reynolds number of 3000}
  \begin{center}
    \includegraphics[scale=0.3]{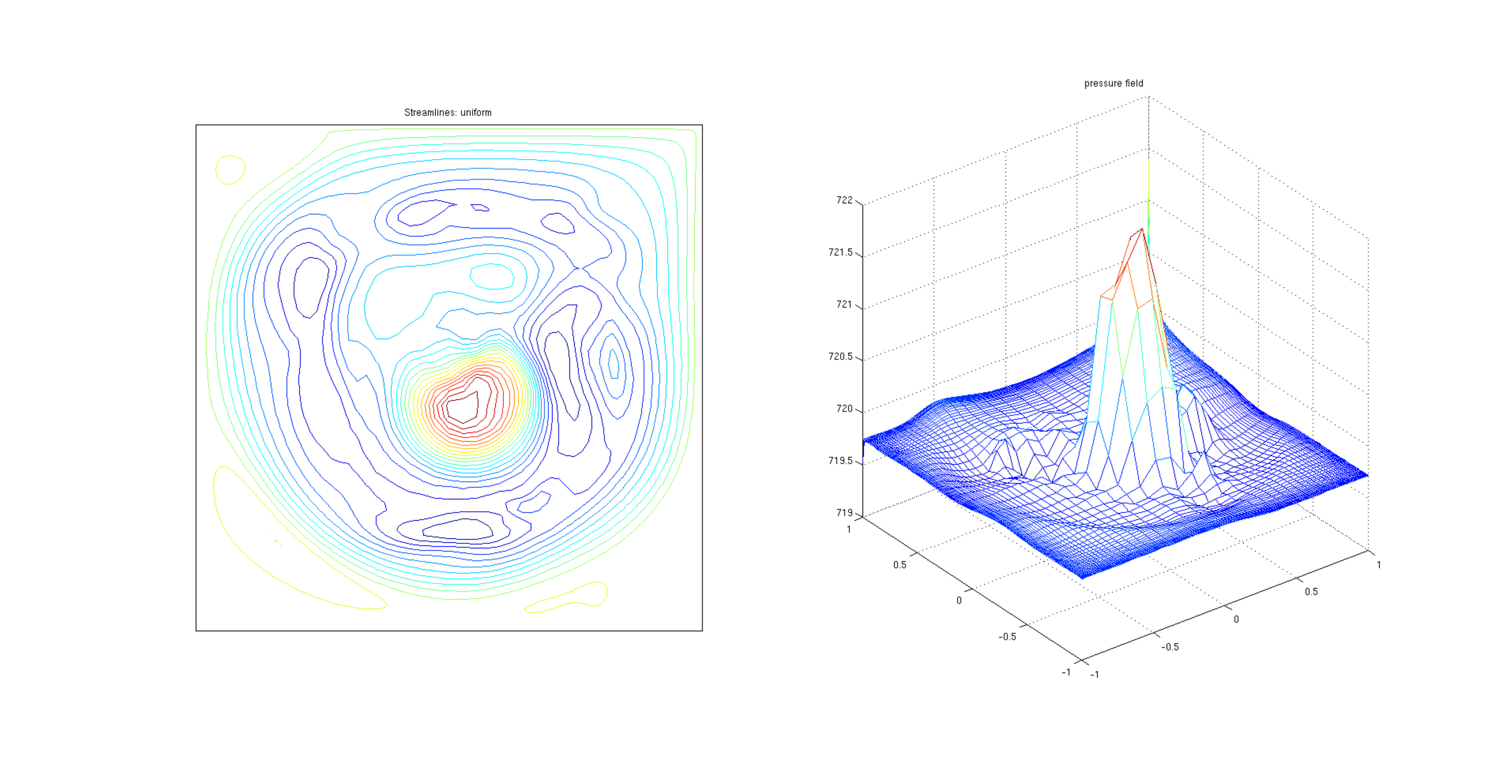} \\
  \end{center}
\end{figure}

In Table (\ref{lidCavityUniform}), we present the results for the
uniform grid. We compare all three important aspects of the methods
proposed namely iteration count, CPU time, and fill factor which is
the ratio of the non non-zeros in the preconditioner and non-zeros in
the original coefficient matrix. We compare new methods with the
Pressure convection diffusion and Least square commutator methods as
implemented in IFISS MATLAB toolbox. In the tables, $nA$ denote the
number of diagonal block in (1,1) block and $nS$ denote the number of
independent Schur complement computations. We use an incomplete LU
with tolerance $10^{-4}$ as an inexact solver for the (1,1) blocks.
As expected the number of iterations for the PCD and LSC methods seem
to be independent of the mesh size $h$. But for the new MSC based
methods there is slight increase in the iteration count. This is
expected for a domain decomposition based method; more the number of
subdomains more the iteration count, but more the
parallelism. Nevertheless, the CPU times of the new methods are better
compared to both PCD and LSC with many independent blocks. For
instance, for the $h=1/32$ we have 5 subdomains for the (1,1) block
and 11 subdomains for the Schur complement, while for $h=1/128$, we
have 12 subdomains for the (1,1) block and 31 independent Schur
complement computation.  On comparing the dependence of iteration
count with increasing Reynolds number, we observe that iteration count
increses mildy until Reynolds number is equal to 1000. We notice that
the system remains in steady state untill Reynolds number of 1000 as
seen in figure (\ref{stream2}). For Reynold number larger than 1000,
we expect the unsteady state with several vortices in the strealine
contours as shown in figure (\ref{stream1}); the velocity contours
suggests that the motion is nearly chaotic.
Comparing with PCD and LSC for
Reynolds number higher that 500 the iteration count for PCD and LCS
are higher that the MSC based preconditioners. The size of each of the
MSC are kept equal and they are equal to $sz$ in the table. The
overlap for the OMSCN and OMSCNR methdods are equal to $sz$. Clearly,
the overlap shows some improvements since, the approximation is better
compared to the non-overlaping case. We have an interesting
observation that for LUM shows same iteration count as MSCN, the
reason is that the $E_{ii}$ and $F_{ii}$ blocks are very sparse and
often zero and thus $G_{ii}-F_{ii}D_{ii}^{-1}E_{ii}\approx
G_{ii}$. However, this is only true for the uniform grid, for
stretched grid in Table (\ref{lidCavityStretched}) MSCN converges
faster compared to LUM method. But all in all for both the grid OMSCN
and OMSCNR seems to have less iteration count and less CPU time
compared to the PCD and LSC methods. For MSCE method the convergence
time was found to be very large compared to all the methods and it
will be investigated in future work.
\begin{sidewaystable}
  \small
  \caption{ \label{lidCavityUniform} Preconditioned GMRES on steady
    Oseen problems, leaky lid driven cavity, (Q2-Q1 FEM, uniform
    grids), for MSCN, LUM, MSCE, OMSCN, OMSCNR, MPCD, and LSC }
  
  \begin{center}
    \begin{tabular}{llllll lll lll lll lll lll ll ll}
      \hline
      Re&  $\frac{1}{h}$ & sz  & nA & nS & sA  & \multicolumn{3}{c}{MSCN} & \multicolumn{3}{c}{LUM} & \multicolumn{3}{c}{MSCE}&\multicolumn{3}{c}{OMSCN} & \multicolumn{3}{c}{OMSCNR} & \multicolumn{2}{c}{MPCD} & \multicolumn{2}{c}{LSC} \\
      \cline{7-9} \cline{10-12} \cline{13-15} \cline{16-18} \cline{19-21} \cline{22-23} \cline{24-25}

      &          &         &    &     &      & its                      & tm                    & ff                                                  & its                        & tm                     & ff  & its & tm & ff  & its & tm & ff & its & tm & ff & its & tm & its & tm \\

      &         &   &    &  &   &                        &                      &                                                  &                          &                       &  &    &    &  &  &  &  &  &  &  &  &  &  &  \\

      & 32 & 38  & 5 & 11 & 400 &  64 & 0.9 & 3.1 & 64 & 0.9 & 3.1 & - & - & - & 52 & 0.7 & 3.1 & 52 & 0.7 & 3.1 & 24 & 1.7 & 13 & 0.9 \\
      10 & 64 & 69 & 10 & 16 & 800 &  101 & 8 & 4.6 & 101 & 8 & 4.6 & - & - & - & 57 & 4.8  & 4.6 & 52  & 4.3 & 4.6 & 25 & 7.8 & 17 & 4.5 \\

      &  128   & 107  & 12  & 31 & 2700   & 141                       & 63                     & 6.9                                                 & 141                         & 63                      & 6.9 & -   & -   & - & 104 & 50 &7.0  & 94 & 44 & 7.0 & 25 & 31 & 21 & 26 \\

      &         &   &   &  &   &                        &                      &                                                  &                          &                       &  &    &    &  &  &  &  &  &  &  &  &  &  &  \\

      & 32 & 38  & 5 & 11 & 400  & 66 & 0.9 & 3.5 & 66 & 0.9 & 3.5 & - & - & - & 53 & 0.7 & 3.5 & 53 & 0.7 & 3.5 & 46 & 2.6 & 31 & 2.2 \\
      100 & 64 & 69 & 10 & 16 & 800 & 103 & 7.0 & 4.1 & 103 & 7.0 & 4.1 & - & - & - & 61 & 4.2 & 4.1 & 61 & 4.2 & 4.1 & 48 & 12.5 & 35 & 9.1 \\

      & 128 & 107 & 12 & 31 & 2700  & 165 & 66 & 6.1 & 165 & 66 & 6.1 & - & - & - & 118 & 46 & 6.2 & 116 & 45 & 6.2 & 49 & 55 & 49 & 56 \\

      &         &   &    &  &   &                        &                      &                                                  &                          &                       &  &    &    &  &  &  &  &  &  &  &  &  &  &  \\

      & 32  & 38  & 5 & 11 & 400  & 80 & 1.2 & 4.2 & 80 & 1.1 & 4.2 & - & -  & -   & 66 & 10.4 & 4.2 & 66 & 10.5 & 4.2 & 120 & 8.8 & 96 & 5.5 \\
      500 & 64  & 69 & 10 & 16 & 800  & 106 & 9.1 & 5.8 & 106 & 9.1 & 5.8 & - & -  & -   & 64 & 5.7 & 5.9 & 67 & 5.9 & 5.9 & 117 & 28.4 & 101 & 26.9 \\
      & 128 & 107  & 12 & 31 & 2700  & 161 & 86 & 8.7 & 161 & 83 & 8.7 & - & - & - & 108 & 58 &8.8  &108  & 58 & 8.8 & 108 & 119 & 100 & 121 \\
      &     &     &       &     &     &    &     &     &    &     &     &   &    &     &     &      &    &     &      &    &     &      &     &      \\

      & 32  & 38  & 5 & 11 & 400  & 96 & 1.4 & 4.1 & 96 & 1.4 & 4.1 & - & -  & -   & 83 & 1.2 & 4.2 & 83 & 1.2 & 4.2 & 179 & 10.5 & 155 & 9.0 \\
      1000 & 64  & 69 & 10 & 16 & 800  & 136 & 11 & 6.1 & 136 & 11 & 6.1 & - & -  & -   & 82 & 7.4 & 6.2 & 91 & 8.1 & 6.2 & 213 & 53.2 & 191 & 52.1 \\
      & 128 & 107  & 12 & 31 & 2700  & 141 & 63 & 6.9 & 141 & 66 & 6.9 & - & - & - & 104 & 49 & 7.0 & 94 & 44 & 7.0 & 198 & 225 & 187 & 231 \\

      &     &     &       &     &     &    &     &     &    &     &     &   &    &     &     &      &    &     &      &    &     &      &     &      \\

      & 32  & 38  & 5 & 11 & 400  & 253 & 5.1 & 4.7 & 253 & 5.1 & 4.7 & - & -  & -   & 232 & 4.6 & 4.8 & 232 & 4.6 & 4.8 & 302 & 18.5 & 300 & 18.6 \\
      3000 & 64  & 69 & 10 & 16 & 800  & 281 & 31 & 6.9 & 281 & 32 & 6.9 & - & -  & -   & 183 & 19 & 7.0 & 183 & 19 & 7.0 & 1127 & 362.3 & 1103 & 381.8 \\
      & 128 & 107  & 12 & 31 & 2700  & 984 & 900 & 11.4 & 984 & 900 & 11.4 & - & - & - & 854 & 705 & 11.4 & 847 & 761 & 11.4 &NC & NC & NC & NC \\

      \hline
    \end{tabular}
  \end{center}
\end{sidewaystable}

\begin{sidewaystable}
  \small

  \caption{ \label{lidCavityStretched} Preconditioned GMRES on steady
    Oseen problems, leaky lid driven cavity, (Q2-Q1 FEM, stretched
    grids), for MSCN, LUM, MSCE, OMSCN, OMSCNR, MPCD, and LSC }

  \begin{center}
    \begin{tabular}{llllll lll lll lll lll lll ll ll}
      \hline
      Re&  $\frac{1}{h}$ & sz  & nA & nS & sA  & \multicolumn{3}{c}{MSCN} & \multicolumn{3}{c}{LUM} & \multicolumn{3}{c}{MSCE}&\multicolumn{3}{c}{OMSCN} & \multicolumn{3}{c}{OMSCNR} & \multicolumn{2}{c}{MPCD} & \multicolumn{2}{c}{LSC} \\
      \cline{7-9} \cline{10-12} \cline{13-15} \cline{16-18} \cline{19-21} \cline{22-23} \cline{24-25}

      &          &         &    &     &      & its                      & tm                    & ff                                                  & its                        & time                     & ff  & its & tm & ff  & its & tm & ff & its & tm & ff & its & tm & its & tm \\

      &         &   &  &  &  &   &                        &                                                                        &                          &                       &  &    &    &  &  &  &  &  &  &  &  &  &  &  \\

      & 32 & 25  & 5 & 17 & 400  & 74 & 0.9 & 3.0 & 74 & 0.9 & 3.0 & - & - & - & 46 & 0.6 & 3.0 & 46 & 0.6 & 3.0 & 23 & 1.3 & 10 & 1.0 \\
      10 & 64 & 101 & 9 & 11 & 900  & 65 & 4.7 & 4.1 & 86 & 6.3 & 4.1 & - & - & - & 49 & 3.7  & 4.1 & 48  & 3.6 & 4.1 & 23 & 8.0 & 28 & 9.2 \\

      &  128   & 238  & 12  & 14 & 2700   & 121                       & 51                     & 6.1                                                 & 113                         & 47                      & 6.1 & -   & -   & - & 79 & 35 &6.2  & 85 & 38 & 6.2 & 23 & 29 & 40 & 49.2 \\

      &         &   &  &  &  &                           &                      &                                                  &                          &                       &  &    &    &  &  &  &  &  &  &  &  &  &  &  \\

      & 32 & 25  & 5 & 17 & 400  & 80 & 1.1 & 3.2 & 80 & 1.1 & 3.2 & - & - & - & 53 & 0.7 & 3.3 & 53 & 0.7 & 3.3 & 47 & 3.3 & 41 & 2.3 \\
      100 & 64 & 101 & 9 & 11 & 900  & 64 & 4.2 & 3.8 & 87 & 5.7 & 3.8 & - & - & - & 49 & 3.3 & 3.9 & 49 & 3.3 & 3.8 & 47 & 13.7 & 62 & 15.9 \\

      & 128 & 238 & 12 & 14 & 2700  & 120 & 46 & 5.5 & 142 & 58 & 5.5 & - & - & - & 99 & 39 &5.6  & 101 & 40 & 5.6 & 48 & 60 & 92 & 117 \\

      &         &   &  &  &  &                           &                      &                                                  &                          &                       &  &    &    &  &  &  &  &  &  &  &  &  &  &  \\

      & 32  & 25  & 5 & 17 & 400  & 107 & 1.6 & 4.2 & 107 & 1.6 & 4.2 & - & -  & -   & 86 & 8.8 & 4.3 & 86 & 8.8 & 4.3 & 160 & 11.1 & 99 & 7.1 \\
      500 & 64  & 101 & 9 & 11 & 900  & 65 & 5.2 & 5.1 & 93 & 7.3 & 6.1 & - & -  & -   & 46 & 4.0 & 5.2 & 48 & 4.1 & 5.2 & 128 & 34.4 & 143 & 37.6 \\
      & 128 & 238  & 12 & 14 & 2700  & 128 & 55 & 6.9 & 170 & 78 & 6.9 & - & - & - & 106 & 47 &7.0  &114  & 51 & 7.0 & 133 & 164 & 193 & 235 \\
      &     &     &   &    &          &    &     &     &    &     &     &   &    &     &     &      &    &     &      &    &     &      &     &      \\

      & 32  & 25  & 5 & 17 & 400  & 120 & 1.9 & 4.2 & 120 & 1.9 & 4.2 & - & -  & -   & 95 & 1.4 & 4.3 & 95 & 1.5 & 4.3 & 251 & 15.2 & 153 & 11.3 \\
      1000 & 64  & 101 & 9 & 11 & 900  & 78 & 6.5 & 5.5 & 111 & 9.3 & 5.5 & - & -  & -   & 59 &5.2  &5.6  & 59 & 5.2 & 5.6 & 248 & 69.3 & 229 & 61.3 \\
      & 128 & 238 & 12 & 14 & 2700  & 133 & 65 & 7.9 & 176 & 88 & 7.9 & - & - & - & 115 & 56 & 8.0 & 119 & 58 & 7.9 & 286 & 364 & 318 & 393 \\

      &     &     &   &    &         &    &     &     &    &     &     &   &    &     &     &      &    &     &      &    &     &      &     &      \\

      & 32  & 25  & 5 & 17 & 400  & 225 & 4.6 & 4.8 & 225 & 4.6 & 4.8 & - & -  & -   & 175 & 3.3 & 4.9 & 175 & 3.4 & 4.9 & 297 & 18.6 & 299 & 19.6 \\
      3000 & 64  & 101 & 9 & 11 & 900  & 240 & 23 & 6.2 & 284 & 29 & 6.2 & - & -  & -   & 216 &21.4  &6.3  & 226 & 22.2 & 6.3 & 940 & 316 & 766 & 240 \\
      & 128 & 238 & 12 & 14 & 2700  & 184 & 101 & 9.3 & 232 & 130 & 9.3 & - & - & - & 163 & 90 & 9.4 & 165 & 91 & 9.4 & 1513 & 2228 & 1232 & 1726 \\
      \hline
    \end{tabular}
  \end{center}
\end{sidewaystable}

\normalfont


\section{Conclusions}
In this work, we proposed a new class of algebraic approximation to
the Schur complement based preconditioners for Navier-Stokes problem,
we observe the superiority of the methods compared to PCD and LSC
method especially for Reynolds number larger than 100. The proposed
methods are highly parallel, and very suitable for modern day
multiprocessor and/or multi-core systems.

As a future work, an implementation of the methods in parallel may be
done and an extensive comparison with other existing preconditioners
for the Navier-Stokes systems will be added. A coarse grid corrections
may be introduced which may improve the convergence of the methods
significantly, however, it may very well be the bottleneck in
achieving overall parallelism in the two level methods.

\section{Acknowledgements}

I am pleased to acknowledge the generous support of Institute Henri
Poincare, Paris and Universit\'e libre de Bruxelles for usual office
facilities, access to library, and encouraging environment that helped
in completing this work.


\end{document}